\documentclass[a4paper, english]{amsart}

\usepackage{color}
\usepackage{amsmath}
\usepackage{booktabs}
\usepackage{amssymb}
\usepackage{wasysym}
\usepackage{graphicx}
\usepackage{verbatim}
\usepackage{mathrsfs}
\usepackage[all]{xy}
\usepackage{enumerate}
\usepackage{hyperref}

\newtheorem{MainThm}{Theorem}

\newtheorem{thm}{Theorem}[section]
\newtheorem{cor}[thm]{Corollary}
\newtheorem{lem}[thm]{Lemma}
\newtheorem{prop}[thm]{Proposition}
\theoremstyle{definition}
\newtheorem{defn}[thm]{Definition}
\newtheorem{example}[thm]{Example}
\theoremstyle{remark}
\newtheorem{rem}[thm]{Remark}
\numberwithin{equation}{section}

\newcommand{\bC}{\mathbb{C}}
\newcommand{\bD}{\mathbb{D}}
\newcommand{\bG}{\mathbb{G}}
\newcommand{\bN}{\mathbb{N}}

\newcommand{\bQ}{\mathbb{Q}}
\newcommand{\bR}{\mathbb{R}}
\newcommand{\bS}{\mathbb{S}}
\newcommand{\bT}{\mathbb{T}}
\newcommand{\bZ}{\mathbb{Z}}

\newcommand{\cB}{\mathcal{B}}
\newcommand{\cE}{\mathcal{E}}
\newcommand{\cF}{\mathcal{F}}
\newcommand{\cG}{\mathcal{G}}
\newcommand{\cO}{\mathcal{O}}

\newcommand{\cT}{\mathcal{T}}

\newcommand{\MT}[2]{\mathbf{MT #1}(#2)}
\newcommand{\MTSO}{\mathbf{MTSO}}

\newcommand\lra{\longrightarrow}

\newcommand\colim{\mathrm{colim \,}}

\newcommand{\hcoker}{/\!\!/}

\newcommand{\hofib}{\operatorname{hofib}}
\newcommand{\hq}{/\hspace{-1.2mm}/}
\newcommand{\delbar}{\overline{\partial}}
\newcommand{\cp}{\mathbb{CP}}
\newcommand{\Sym}{\operatorname{Sym}}

\newcommand{\mclass}{\kappa}
\newcommand{\nclass}{\nu}

\newcommand{\CircNum}[1]{\ooalign{\hfil\raise .00ex\hbox{\scriptsize #1}\hfil\crcr\mathhexbox20D}}

\newcommand{\X}{\mathbf{X}}
\newcommand{\map}{\mathrm{map}}
\newcommand{\CPinf}{B\bC^\times}
\newcommand{\loopinf}{\Omega^{\infty}}

\newcommand{\Top}{\mathbf{Top}}
\newcommand{\Hol}{\mathrm{Hol}}
\newcommand{\Hom}{\mathrm{Hom}}
\newcommand{\Pic}{\mathrm{Pic}}
\newcommand{\HolStack}{\mathsf{Hol}}
\newcommand{\CxStack}{\mathsf{Cx}}
\newcommand{\ModStack}{\mathsf{M}}
\newcommand{\PicStack}{\mathsf{Pic}}
\newcommand{\qHolStack}{\mathsf{qHol}}
\newcommand{\qHolstack}{\mathsf{qHol}}
\newcommand{\mtmap}{\alpha}
\newcommand{\sect}{\Gamma}
\newcommand{\Gl}{\operatorname{Gl}}
\newcommand{\pgl}{\mathbb{P}\Gl}

\newcommand{\ch}{\operatorname{ch}}
\newcommand{\ind}{\operatorname{Ind}}
\newcommand{\Td}{\operatorname{Td}}
\newcommand{\pr}{\mathrm{pr}}

\newcommand{\modspace}[1]{\mathcal{S}_{#1}}
\newcommand{\moddspace}[2]{\mathcal{S}_{#1}^{#2}}

\newcommand{\Diff}[1]{\mathscr{D}_{#1}}
\newcommand{\Difff}[2]{\mathscr{D}_{#1}^{#2}}
\newcommand{\exmcg}{\widetilde{\Gamma}}

\renewcommand{\top}{\mathrm{top}}
\newcommand{\hol}{\mathrm{hol}}

\mathchardef\ordinarycolon\mathcode`\:
\mathcode`\:=\string"8000
\begingroup \catcode`\:=\active
  \gdef:{\mathrel{\mathop\ordinarycolon}}
\endgroup

\title[Universal Picard varieties and the extended mapping class group]{Stable cohomology of the universal Picard varieties and the extended mapping class group}
\author{Johannes Ebert}
\email{johannes.ebert@uni-muenster.de}
\address{Mathematisches Institut der Westf\"alischen Wilhelms-Universit\"at M\"unster\\
Einstein\-stra{\ss}e 62\\
DE-48149 M\"unster\\
Germany}

\author{Oscar Randal-Williams}
\email{o.randal-williams@math.ku.dk}
\address{Institut for Matematiske Fag\\
Universitetsparken 5\\
DK-2100 K{\o}benhavn {\O}\\
Denmark}
\date{\today}

\begin{document}

\begin{abstract}
We study the moduli spaces which classify smooth surfaces along with a complex line bundle. There are homological stability and Madsen-Weiss type results for
these spaces (mostly due to Cohen and Madsen), and we discuss the
cohomological calculations which may be deduced from them. We then relate these spaces to (a generalisation of) Kawazumi's extended mapping class groups, and hence deduce cohomological information about these. 

Finally, we relate these results to complex algebraic geometry. We construct a holomorphic stack classifying families of Riemann surfaces equipped with a fibrewise holomorphic line bundle, which is a gerbe over the universal Picard variety, and compute its holomorphic Picard group.
\end{abstract}
\maketitle

\section{Introduction}

Let $\Sigma_{g,r}^{n}$ be a connected oriented smooth surface of genus $g$ with $r$ boundary components and $P:=\{p_1,\ldots,p_n \}$ be $n$ distinct marked points in the interior. Let $\Difff{g,r}{n}$ the group of orientation preserving diffeomorphisms of $\Sigma_{g,r}^{n}$ which fix the set $\partial \Sigma_{g,r}^n \cup P$ pointwise and let $\Gamma_{g,r}^{n}:=\pi_0 \Difff{g,r}{n}$ be the mapping class group.
We assume that $g \geq 2$, so that by a theorem of Earle and Eells \cite{EE} $\Difff{g,r}{n}$ has contractible components and the natural map on classifying spaces $B \Difff{g,r}{n} \to B \Gamma_{g,r}^{n}$ is a homotopy equivalence.

Let $\map_{\partial} (\Sigma_{g,r}^{n}, \CPinf)$ denote the space of continuous maps $\Sigma_{g,r}^{n} \to \CPinf$ which send $\partial \Sigma_{g,r}^n \cup P$ to the basepoint of $\CPinf$. We write $\map_{\partial} (\Sigma_{g,r}^{n}, \CPinf) (k)$ for the path component consisting of those maps of degree $k$. The principal object of study in this paper is the moduli space of surfaces of genus $g$ with $r$ boundary components, $n$ marked points and a complex line bundle, originally introduced by Cohen and Madsen \cite{CM}:
\begin{equation*}
\moddspace{g,r}{n}:= \map_{\partial}(\Sigma_{g,r}^{n}, \CPinf) \times_{\Difff{g,r}{n}} E \Difff{g,r}{n}.
\end{equation*}
Being a Borel construction, this space fits into a \emph{defining fibration sequence}
\begin{equation}\label{eq:DefiningFibration}
\map_{\partial}(\Sigma_{g,r}^{n}, \CPinf) \lra \moddspace{g,r}{n} \lra B\Difff{g,r}{n}.
\end{equation}
The space $\moddspace{g,r}{n}$ carries a tautological surface bundle
\begin{equation*}
\mathcal{E}_{g,r}^{n}:= (\Sigma_{g,r}^n \times \map_{\partial}(\Sigma_{g,r}^{n}, \CPinf)) \times_{\Difff{g,r}{n}} E \Difff{g,r}{n},
\end{equation*}
and we denote the bundle projection by $\pi : \mathcal{E}_{g,r}^{n} \to \moddspace{g,r}{n}$, the vertical tangent bundle by $T_v$, and the $n$ canonical sections by $s_1, ..., s_n$. The evaluation map $\Sigma_{g,r}^n \times \map_{\partial}(\Sigma_{g,r}^{n}, \CPinf) \to \CPinf$ is $\Difff{g,r}{n}$-invariant, so gives a map $L:\mathcal{E}_{g,r}^{n} \to \CPinf$. We use the letter $L$ also for the complex line bundle this map induces, and note that it is trivialised over the boundary and marked points in each fibre.

The space $\moddspace{g,r}{n}$ splits into path components $\moddspace{g,r}{n}(k)$ corresponding to the maps of degree $k$, and the total space $\mathcal{E}_{g,r}^{n}$ has an analogous decomposition. The term ``moduli space'' is justified by the following observation. 

\begin{prop}\label{classifying-space}
Let $X$ be any paracompact Hausdorff space. Homotopy classes of maps $X \to \moddspace{g,r}{n} (k)$ are in bijection with isomorphism classes of tuples $(E, \pi,s ,L,\eta)$, where $\pi:E \to X$ is an oriented smooth surface bundle of genus $g$ with $r$ boundary components and the boundary bundle is trivialised; $s=(s_1,\ldots,s_n)$ are pairwise disjoint cross-sections in the interior; $L \to E$ is a line bundle; $\eta$ is a trivialisation of $L$ over $\partial E \cup s_1 (X) \cup \ldots \cup s_n (X)$, such that with respect to this trivialisation, $L$ has fibrewise degree $k$.
\end{prop}

\begin{defn}\label{def:extendedMCG}
For each $k$ we choose a degree $k$ map $L_0 \in \map_\partial(\Sigma_{g,r}^n, \CPinf)(k)$. The \emph{extended mapping class group} $\exmcg_{g,r}^{n}(k)$ is the fundamental group of the space $\moddspace{g,r}{n}(k)$, based at $L_0$,
$$\exmcg_{g,r}^{n} (k):= \pi_1 (\moddspace{g,r}{n}(k), L_0).$$
\end{defn}

When $k=0$ we make the canonical choice of $L_0$, which is the constant map to the basepoint. In this case the extended mapping class group also has an algebraic description. If we let $H_{g,r}^{n}$ denote the $\Gamma_{g,r}^{n}$-module $H^1 (\Sigma_{g,r}^{n}, \partial \Sigma_{g,r}^{n} \cup P;\bZ)$, there is an isomorphism 
$$\exmcg_{g,r}^{n}(0)\cong H_{g,r}^{n} \rtimes \Gamma_{g,r}^{n}.$$
If $(r,n,k) = (1,0,0)$ this definition agrees with that of Kawazumi \cite{Kawazumi98}, however we warn the reader that our definition is different from Kawazumi's for other values of $(r,n)$.

For both $\moddspace{g,r}{n}(k)$ and $\exmcg_{g,r}^n(k)$ we adopt the convention of omitting $r$ or $n$ from the notation when they are zero.

\subsection{Topological results}

The first part of this paper (Sections \ref{sec:moduli-spaces} and \ref{sec:cohomology-exmcg}) is devoted to the study of the cohomology of the group $\exmcg_{g,r}(k)$. We will see in Proposition \ref{homotopy-type-modspace} that as long as $r+n >0$ the space $\moddspace{g,r}{n}(k)$ is aspherical, so there is a homotopy equivalence 
\begin{equation*}
 B \exmcg_{g,r}^{n} (k) \simeq \moddspace{g,r}{n}(k),
\end{equation*}
and the cohomology of the space $\moddspace{g,r}{n}(k)$ translates to the group cohomology of $\exmcg_{g,r}^{n} (k)$. The relation between $\modspace{g}(k)$ and $B \exmcg_g (k)$ turns out to be much more subtle.

Due to the moduli-theoretic interpretation of $\moddspace{g,r}{n} (k)$, the methods of the homotopy theory of moduli spaces apply and yield a great amount of information about these spaces, and hence about $B \exmcg_{g,r}^{n}(k)$ for $r+n>0$. Comparison with $\moddspace{g}{}(k)$ will also be essential for the study of $B \exmcg_g (k)$.

\begin{MainThm}\label{modspaceproperties}
\leavevmode
\begin{enumerate}[(i)]
\item \label{it:independence} (Independence of $k$) There are homotopy equivalences
\begin{align*}
\moddspace{g,r}{n}(k) &\simeq \moddspace{g,r}{n}(k+1)  \; \text{for}\; r>0\\
\moddspace{g}{1} (k) &\simeq \moddspace{g}{1} (k+1)\\
\modspace{g}(k) &\simeq \modspace{g}(k+2g-2).
\end{align*}
\item \label{it:stability} (Harer type stability) Consider the ``stabilisation maps"
\begin{align*}
\alpha(g)&: \moddspace{g,r+1}{n}(k) \lra \moddspace{g+1,r}{n}(k)\\
\beta(g) &:  \moddspace{g,r}{n}(k) \lra \moddspace{g,r+1}{n}(k)\\
\gamma(g)&: \moddspace{g,r}{n}(k) \lra \moddspace{g,r-1}{n}(k)
\end{align*}
given by gluing in a pair of pants along the legs, along the waist, or gluing in a disc (and extending the map to $\CPinf$ trivially over the glued-in part in each case). Then the induced maps in homology $H_* (\alpha(g))$, $H_* (\beta(g))$, $H_* (\gamma(g))$ are isomorphisms if $3* \leq 2g-3$, $3* \leq 2g-1$,  $3* \leq 2g-1$, respectively.
\item \label{it:splitting} (Splitting) The map $\moddspace{g,r}{n}(k) \to \modspace{g,r}(k) \times (\CPinf)^n$ given by the tangent spaces at the marked points is a homology isomorphism in degrees $3* \leq 2g-1$.
\item \label{it:MW} (Madsen--Weiss type theorem) Let $\MTSO(2)$ be the Madsen--Tillmann--Weiss spectrum. There is a map $\mtmap_{g,r} : \modspace{g,r}(k) \to \loopinf_0 (\MTSO(2) \wedge \CPinf_+)$, which is a homology equivalence in degrees $3* \leq 2g-3$.
\end{enumerate}
\end{MainThm}
We shall say ``stable range'' for the range of degrees in which a given homology group is independent of $g$. 

Most of this theorem has been proved elsewhere. A complete proof of homological stability as formulated in (\ref{it:stability}) was given by the second named author \cite{R-WResolution}, but the result has a longer prehistory. For $r>0$, it was first proved by Cohen--Madsen \cite{CM}, but with the weaker stable range $2* \apprle g$. The range of stability was improved by Boldsen \cite{Boldsen}, again as long as $r>0$. The proof of stability for closing the last boundary in \cite{CM} is incorrect and has been partially repaired by those authors in \cite{CM2}, as long as rational coefficients are employed.

In the references we have cited only the case $n=0$ is covered, but in Section \ref{sec:MarkedPoints} we show how to derive the case $n>0$. The homological splitting of (\ref{it:splitting}) is a straightforward adaption of an argument by B\"odigheimer and Tillmann \cite{BT}, and we briefly sketch it in Section \ref{sec:MarkedPoints}. Part (\ref{it:MW}) is a generalisation of the Madsen--Weiss theorem \cite{MW} due to Cohen--Madsen \cite{CM} (though one needs the result of \cite{R-WResolution} to deal with the case $r=0$), and we recall this result in Section \ref{sec:MadsenWeiss}. We prove (\ref{it:independence}) in Section \ref{sec:periodicity}. 
\begin{rem}
Parts (\ref{it:stability}) and (\ref{it:splitting}) show that $H_* (\moddspace{g}{n}(k) )\cong H_* (\moddspace{g}{n}(k+1) )$ in the stable range. However, in general $\modspace{g} (k) \not \simeq \modspace{g}(k+1)$. This can be deduced from Theorem \ref{thm:LowDimCohomology}, which computes the group cohomology of the fundamental groups of these spaces and we see they differ in general.
\end{rem}

The space $\loopinf_0 (\MTSO(2) \wedge \CPinf_+)$ can be understood by the usual methods of algebraic topology. We defer its definition to Section \ref{sec:MadsenWeiss}, where we will also construct certain cohomology classes
$$\mclass_{i,j} \in H^{2i+2j}(\Omega^\infty_0 (\MTSO(2) \wedge \CPinf_+);\bZ).$$
Under the map $\alpha_{g,r} : \modspace{g,r}(k) \to \loopinf_0 (\MTSO(2) \wedge \CPinf_+)$ the class $\mclass_{i,j}$ pulls back to
$$\pi_{!} (e(T_v)^{i+1} c_{1}(L)^{j}) \in H^{2i+2j}(\moddspace{g,r}{n}(k);\bZ),$$
where $\pi : \mathcal{E}_{g,r}(k) \to \modspace{g,r}(k)$ is the tautological surface bundle, $T_v$ is the vertical tangent bundle and $L$ is the complex line bundle we described earlier. For simplicity we denote this class by $\mclass_{i,j}$ too. The class $\mclass_{i,0}$ may be defined on $B \Diff{g,r}$, and here it coincides with the Mumford--Morita--Miller  class usually denoted $\kappa_i$.

\begin{MainThm}\label{thm:HolRationalCohomology}
There is an isomorphism
$$H^*(\loopinf_0 (\MTSO(2) \wedge \CPinf_+);\bQ) \cong  \bQ[\mclass_{i,j} \,\, | \,\, i+j > 0, j\geq 0, i \geq -1].$$
\end{MainThm}

Therefore when $r > 0$, for degrees $3* \leq 2g-3$ there are isomorphisms
\begin{equation}\label{modspaceRationalcohomology}
H^ * (B \exmcg_{g,r} (k); \bQ) \cong H^*(\modspace{g,r}(k);\bQ) \cong  \bQ[\mclass_{i,j} \,\, | \,\, i+j > 0, j\geq 0, i \geq -1].
\end{equation}

We are also able to make low-dimensional integral calculations. The \emph{Hodge class} $\lambda \in H^2(\modspace{g};\bZ)$ may be defined as $c_1(\pi_!^K(T_v^*))$, the first Chern class of the push-forward in complex $K$-theory of the vertical cotangent bundle of the tautological surface bundle, and it is known to satisfy the relation $12\lambda = \kappa_1$. Just as for the $\mclass_{i,j}$, the class $\lambda$ is induced from a class on $\loopinf (\MTSO(2) \wedge \CPinf_+)$ which we also denote $\lambda$.

\begin{MainThm}\label{thm:LowDimCohomologyHol}
The groups
$$H^1(\loopinf_{0} (\MTSO(2) \wedge \CPinf_+);\bZ) \quad \text{and} \quad H^3(\loopinf_{0} (\MTSO(2) \wedge \CPinf_+);\bZ)$$
are trivial. The group $H^2(\loopinf_{0} (\MTSO(2) \wedge \CPinf_+);\bZ)$ is free abelian of rank three, with free basis the Hodge class $\lambda$, $\mclass_{0,1}$ and $\zeta := \tfrac{1}{2}(\mclass_{0,1} - \mclass_{-1,2})$.
\end{MainThm}

As long as $g \geq 6$, this gives a description of the integral cohomology of $\modspace{g,r}(k)$ up to degree 3. 


 We now turn to the case $r=n=0$. The main obstacle in studying the cohomology of $B \exmcg_g (k)$ is that the stabilisation map $B \exmcg_{g,1}(k) \to B \exmcg_g (k)$ is not a homology isomorphism in a range of degrees. However the failure of homological stability is entirely down to the failure of stability on the third cohomology. Before we state the key result, let us introduce some language.

\begin{defn}
A \emph{fibre sequence} consist of a connected pointed space $(B, b)$, maps $F \stackrel{f}{\to} E \stackrel{p}{\to} B$ and a nullhomotopy $H : p \circ f \simeq c_b$ from the composition to the constant map to $b$. This determines a map $F \to \hofib_b(p)$ which we require to be a weak equivalence. We will often drop the nullhomotopy from the notation.
\end{defn}

Define
$$\Xi:\CPinf \overset{-\otimes L_0}\lra \map(\Sigma_g, \CPinf)(k) \lra \modspace{g}(k),$$
which is a map that classifies the data $(\CPinf \times \Sigma_g, \pr_1, \pr_{1}^{*} \gamma \otimes \pr_{2}^{*} L_0)$, where $\gamma$ denotes the universal line bundle over $\CPinf$. 

\begin{MainThm}\label{thm:ExtendedMCGStability}
\leavevmode
\begin{enumerate}[(i)]
\item For $g \geq 3$, there exist fibre sequences
\begin{eqnarray*}\label{fibration-sequence}
\CPinf \stackrel{\Xi}{\lra} \modspace{g}(k) \stackrel{\Pi}{\lra} B \exmcg_g (k)\\
\modspace{g}(k) \stackrel{\Pi}{\lra} B \exmcg_g (k) \stackrel{\theta}{\lra} K(\bZ,3)
\end{eqnarray*}
(i.e.\ there exist nullhomotopies of $\Pi \circ \Xi$ and $\theta \circ \Pi$ yielding fibre sequences). 
\item There is a map from $B\exmcg_{g,1}(k)$ to the homotopy fibre of the map $\theta$ which induces an isomorphism in integral homology in degrees $3* \leq 2g-1$.
\end{enumerate}
\end{MainThm}






It turns out that the cohomology of $\exmcg_g(k)$ behaves in a fairly systematic way in a certain range of degrees, but the systematic way in which it behaves in turn depends on $k$.
The following are our main results on the cohomology of $\exmcg_g(k)$.

\begin{MainThm}\label{thm:LowDimCohomology}
Suppose $g \geq 6$. The group $H^1(B \exmcg_g(k);\bZ)$ is trivial, and the group $H^2(B \exmcg_g(k);\bZ)$ is free abelian of rank two and injects into $H^2(\modspace{g}(k);\bZ)$. A free basis for it may be taken to be the Hodge class $\lambda$ and an element $\eta$ that maps to
$$\frac{1}{\gcd(2g-2, g+k-1)}\left ( k \mclass_{0,1} + (g-1)\mclass_{-1,2} \right)$$
in $H^2(\modspace{g}(k);\bZ)$. The group $H^3(B \exmcg_g(k);\bZ)$ is $\bZ/\gcd(2g-2, 1-g-k)$ generated by the class $\theta \in H^3 (B \exmcg_g (k);\bZ) = [B \exmcg_g (k), K(\bZ,3)]$ of Theorem \ref{thm:ExtendedMCGStability}.
\end{MainThm}

We are also able to give a complete description of the rational cohomology algebra in the stable range, where the dependence on $g$ and $k$ is less vividly seen.

\begin{MainThm}\label{thm:PicRationalCohomology}
The rational cohomology $H^*(B \exmcg_g(k);\bQ)$ injects into $H^*(\modspace{g}(k);\bQ)$ in all degrees.
There exist classes $\nclass_{i,j} \in H^{2i+2j}(B \exmcg_g (k); \bQ)$ defined for $i\geq -1$ and $j \geq 0$, such that the algebra homomorphism
$$\bQ[\nclass_{i,j} \mid i\geq -1, j \geq 0, i+j >0; (i,j)\neq (0,1)] \lra H^{*}(B \exmcg_{g}(k); \bQ)$$
is an isomorphism in degrees $3* \leq 2g-3$. Moreover, the image of $\nclass_{i,j}$ in $H^* (\modspace{g}(k); \bQ)$ is equal to $\mclass_{i,j}$ modulo the ideal $(\mclass_{0,1})$, and $\nu_{i,0} = \kappa_i$.
\end{MainThm}

\begin{rem}
Using different methods, Kawazumi \cite{Kawazumi98, KawazumiInfinitesimal, KawazumiGoD} obtained the isomorphism (\ref{modspaceRationalcohomology}) when $r=1$ and $k=0$. Looijenga \cite{Looijenga} was able to determine the structure of $H^{*}(B \exmcg_g (k);\bQ)$ as a module over $H^* (B \Gamma_g; \bQ)$, but not the algebra structure.
\end{rem}

Next, we study the group extension
\begin{equation}\label{second-fibration}
0 \lra H_{g,r}^{n} \lra \exmcg_{g,r}^{n} (k)\lra \Gamma_{g,r}^{n} \lra 1
\end{equation}
obtained by taking the fundamental groups of the spaces in the defining fibration (\ref{eq:DefiningFibration}). The $\Gamma_{g,r}^{n}$-modules $H_{g,r}^n = H^1 (\Sigma_{g,r}^{n}, \partial \Sigma_{g,r}^{n} \cup P;\bZ)$ form a \emph{coefficient system} in the sense of \cite{CM} or \cite{Boldsen}. We denote this coefficient system simply by $H$ and its rationalisation by $H_{\bQ}$, thereby explaining notation such as $H^* (\Gamma_{g,r}^ {n};H)$.

\begin{MainThm}\label{thm:DegenerationSpectralSequence}
The rational Leray--Hochschild--Serre spectral sequence of the extension (\ref{second-fibration}) degenerates at $E_2$ in the stable range. If $r+n\leq 1$ it degenerates in all degrees. The associated graded algebra in the stable range is 
\begin{equation*}
\bigoplus_{p, q} H^p(\Gamma_{g,r}^n ; \wedge^q H_\bQ) \cong \bQ[e_1, ,..., e_n, \kappa_1, \kappa_2, \ldots] \otimes \bQ[x_{i,j} \,\, | \,\, i+j >0, j>0, i\geq -1],
\end{equation*}
if $r+n>0$, and
\begin{equation*}
\bigoplus_{p, q} H^p(\Gamma_g ; \wedge^q H_\bQ) \cong \bQ[\kappa_1, \kappa_2, \ldots] \otimes \bQ[x_{i,j} \,\, | \,\, i+j > 0, j>0, i \geq -1, (i,j) \neq (0,1)],
\end{equation*}
if $r+n=0$. In both cases, $e_i$ has bidegree $(p,q)=(2,0)$, $\kappa_i$ has bidegree $(p,q) = (2i,0)$ and $x_{i,j}$ has bidegree $(p,q) = (2i+j,j)$.
\end{MainThm}

The proof of this theorem is via the work of Kawazumi \cite{Kawazumi98, KawazumiInfinitesimal, KawazumiGoD}, who has studied in depth the rational cohomology of $\exmcg_{g,1}$. He has defined certain classes $m_{i,j} \in H^{2i+j-2}(\Gamma_{g,1}, \wedge^j H)$ which are permanent cycles in the spectral sequence
$$E_2^{p,q} = H^p(\Gamma_{g,1}; \wedge^q H_\bQ) \Longrightarrow H^{p+q}(\exmcg_{g,1};\bQ)$$
and detect classes $\widetilde{m}_{i,j} \in H^{2i+2j-2}(\exmcg_{g,1};\bZ)$ he has also defined. He then shows that the $m_{i,j}$ generate $E_2^{*,*}$ in the stable range, so the spectral sequence collapses. If $r=1$ and $n=0$, Theorem \ref{thm:DegenerationSpectralSequence} follows immediately from the more technical statement that Kawazumi's $\widetilde{m}_{i,j}$ agree with our $\kappa_{i-1, j}$. This is somewhat surprising, as the $\widetilde{m}_{i,j}$ are defined using an explicit group cocycle which depends crucially on the surface having boundary, whereas $\kappa_{i-1, j}$ may be defined for closed surfaces. The other cases with $r+n>0$ are derived using homological stability.

\subsection{Relation to complex geometry}

The homotopy types $B \Diff{g}$, $ \modspace{g}(k)$ and $B \exmcg_{g}(k)$ bear a close relationship to the geometry of Riemann surfaces, and the preceding purely topological results can be used to obtain more geometrical results. In order to formulate these results, we have to employ the language of stacks. We adopt the convention of writing stacks as $\mathsf{Abc}$ and their associated homotopy types as $\mathrm{Abc}$. 

Let $\ModStack_g$ be the moduli stack of genus $g$ Riemann surfaces, defined on the site $\Top$ of topological spaces (with the ordinary topology of open covers). Elements of $\ModStack_g (X)$ or, equivalently, maps $X \to \ModStack_g$, correspond to families of genus $g$ Riemann surfaces $\pi:E \to X$. It is well-known, by Teichm\"uller theory, that $\ModStack_g \cong \cT_g \hcoker \Gamma_g$, the orbifold quotient of Teichm\"uller space by $\Gamma_g$. Moreover, the homotopy type of $\ModStack_g$ is just $B \Diff{g}$.

Denote by $\HolStack_g$ the stack over $\Top$ which classifies families of Riemann surfaces of genus $g$, $E \to X$, equipped with a fibrewise holomorphic line bundle $L \to E$. It splits into components $\HolStack_g^k$, where $k$ is the fibrewise degree of $L$. 
Recall that for an individual Riemann surface $S$, the Picard variety $\Pic(S)$ is the complex Lie group of isomorphism classes of holomorphic line bundles on $S$. It splits into components $\Pic^k (S)$ parametrising isomorphism classes of degree $k$ line bundles. There is a (noncanonical) isomorphism $\Pic(S) \cong \bZ \times \Pic^0 (S)$ and $\Pic^0 (S)$ is a complex $g$-dimensional torus. 
If $\pi:E \to X$ is a family of genus $g$ Riemann surfaces, there is a family $\Pic (E/X) \to X$ of complex manifolds whose fibre over $x \in X$ is just the Picard variety $\Pic (E_x)$ of $E_x := \pi^{-1} (x)$.

Denote by $\PicStack_g$ the stack over $\Top$ which classifies families of Riemann surfaces of genus $g$, $E \to X$, equipped with a section $s : X \to \Pic(E/X)$ of the associated bundle of Picard varieties. Again, $\PicStack_g$ splits into components $\PicStack_g^k$, indexed by the fibrewise degree of the line bundle on the total space. Both stacks have a natural map to $\ModStack_g$ that just remembers the underlying Riemann surface. Furthermore, there is a map $\Phi_g^k: \HolStack_g^k \to \PicStack_g^k$ that takes a complex line bundle to its isomorphism class. 

In Theorem \ref{thm:Gerbe} we show that $\Phi_g^k$ is a gerbe with band $\bC^{\times}$. As a consequence, after taking homotopy types we get a fibre sequence
\begin{equation*}
\CPinf \lra \Hol_g^k \overset{\phi_g^k}\lra \Pic_g^k.
\end{equation*}
In Theorem \ref{thm:HolCxStacks} we show that there are equivalences $\Hol_g^k \simeq \modspace{g}(k)$ and $\Pic_g^k \simeq B\exmcg_g(k)$ under which this fibre sequences corresponds to that of Theorem \ref{thm:ExtendedMCGStability}. Thus the cohomological results of the last section also give information about the cohomology of $\HolStack_g^k$ and $\PicStack_g^k$, and of the effect on cohomology of the map $\Phi_g^k$.

Another property of the stacks $\HolStack_g^k$ and $\PicStack_g^k$ is that they are local quotient stacks (this is the class of stacks which are closest to spaces in the sense that the homotopy type of the stacks reflects geometric properties such as the classification of line bundles and gerbes). Therefore, we are able to prove an analogue of a result of Mestrano--Ramanan on the existence of Poincar\'e line bundles.

It is a well-known---but fairly deep---result that the moduli stack $\ModStack_g$ is even a holomorphic stack. This will easily imply that $\HolStack_g^k$ and $\PicStack_g^k$ are holomorphic (local quotient) stacks, so that concepts from holomorphic geometry translate. In particular, the notion of a holomorphic line bundle on these stacks is meaningful. We show that for both these stacks, the holomorphic Picard group is isomorphic to the topological one. The resulting computation of the Picard group of $\PicStack_g^k$ is not new: it is due to Kouvidakis \cite{Kouvidakis} and we recover his result, although we give a different description of the generators. We also discuss the relationship between our results and related results in algebraic geometry, in particular those of Melo and Viviani \cite{MeloViviani}.

\subsection{Acknowledgements}
Both authors are grateful to the MFO Oberwolfach and the Mathematics department at Stanford University for their support with research visits while this work was being carried out. The second author is also grateful to M.\ Melo and F.\ Viviani for explaining their work. O.\ Randal-Williams was supported by an EPSRC Studentship, the EPSRC PhD Plus scheme, ERC Advanced Grant No.\ 228082, and the Danish National Research Foundation through the Centre for Symmetry and Deformation. 

\section{The spaces of surfaces with line bundles}\label{sec:moduli-spaces}

\subsection{Homotopy type}


Recall from the introduction the definition of $\moddspace{g,r}{n}$, the tautological surface bundle $\pi : \mathcal{E}_{g,r}^{n} \to \moddspace{g,r}{n}$, its sections $s_1, ..., s_n$, the vertical tangent bundle $T_v$ and the complex line bundle $L \to \mathcal{E}_{g,r}^{n}$ trivialised along the boundary and at the marked points. 

\begin{prop}\label{homotopy-type-modspace}
Let $g \geq 2$.
\begin{enumerate}[(i)]
\item If $r+n>0$, then $\moddspace{g,r}{n}(k)$ is aspherical, and so $B \exmcg_{g,r}^{n}(k) \simeq \moddspace{g,r}{n}(k)$.
\item If $r=n=0$, then $\pi_i (\modspace{g}(k))=0$ for $i \neq 1,2$ and $\pi_2 (\modspace{g}(k))=\bZ$.
\end{enumerate}
\end{prop}
\begin{proof}
We look at the long exact sequence of the defining fibration sequence 
\begin{equation*}
\map_{\partial}(\Sigma_{g,r}^{n},\CPinf) \lra \moddspace{g,r}{n} \lra B \Difff{g,r}{n}.
\end{equation*}
If $r+n>0$, then $\map_{\partial} (\Sigma_{g,r}^{n}, \CPinf)(k)$ is an Eilenberg--MacLane space of type $B H_{g,r}^{n}$. This is immediate and proves the first part. 
Because $\CPinf$ is an abelian topological group so is $\map(\Sigma_g,\CPinf)$, and hence $\map(\Sigma_g, \CPinf)(k)$ is a product of Eilenberg--MacLane spaces, namely $BH_g \times \CPinf$ (however, this product decomposition is not natural, which is responsible for many of the subtleties of that case). The second part follows.
\end{proof}

\subsection{Independence of the degree}\label{sec:periodicity}

The homotopy equivalences of part (\ref{it:independence}) of Theorem \ref{modspaceproperties} are given by tensoring the line bundle with a fixed line bundle, provided by Proposition \ref{existencelinebundle} below.

\begin{prop}\label{existencelinebundle}
Let $E_{g,r}^{n} \to B  \Difff{g,r}{n}$ be the universal surface bundle, with sections $s_1,\ldots s_n$. Then there is a line bundle $L \to E_{g,r}^{n}$ with fibrewise degree $k$ and such that both $L|_{\partial E_{g,r}^{n}}$ and $s_{i}^{\ast}L$ are trivial in each of the following cases:
\begin{enumerate}[(i)]
\item $r>0$, $n$ arbitrary and $k=1$.
\item $r=0$, $n=k=1$.
\item $r=0$, $n=0$, $k=2-2g$. 
\end{enumerate}
\end{prop}

This proposition can be restated by saying that the map $\moddspace{g,r}{n}(k) \to B \Difff{g,r}{n}$ has a section in the listed cases.

\begin{proof}
The third case is clear, since we can take the vertical tangent bundle $T_v$ (the condition $n=0$ is necessary as we might be unable to trivialise $T_v$ over the cross-sections in the case $n>0$).

In the first case, the boundary bundle is trivialised and so the bundle may be assumed to contain a trivialised collar. Inside the collar, there exists an embedded $B \Gamma_{g,r}^{n} \times \bD^2 \to E_{g,r}^{n}$ over $B \Gamma_{g,r}^{n}$ which does not meet the boundary or the marked points. The standard collapse construction defines a map $h:E_{g,r}^{n} \to B \Gamma_{g,r}^{n} \times \bD^2/\partial \bD^2$ of fibrewise degree $1$. Let $H \to \bD^2 / \partial \bD^2$ be the Hopf bundle and pick a trivialisation over the point $\partial \bD^2 $ and one over the origin of $\bD^2$. The bundle $L:= h^{*} H$, together with the induced trivialisations, does the job in the first case.

For second case, let $U \to B \Gamma_{g}^{1}$ be the closed unit disc bundle in the line bundle $s^{*} T$, where $T$ is the vertical tangent bundle, and $\partial U$ be the unit sphere bundle. There exists a tubular neighborhood, i.e.\ a fibrewise embedding $U \to E_{g}^{1}$ that maps the zero section to $s$. The collapse construction, performed fibrewise, yields a fibrewise map $E_{g}^{1} \to U/\partial U$ of fibrewise degree $1$. Consider the rotation action of $U(1)$ on $\bS^2$ (considered as the Riemann sphere), fixing the points $0, \infty$, and the $\bS^2$-bundle $p:F:=EU(1) \times_{U(1)} \bS^2 \to BU(1)$. There is a bundle map $U/\partial U \to F$, covering the map $\epsilon$ from (\ref{tangentspacemap}). The composition with the collapse is a fibre-preserving map $h: E_{g}^{1} \to F$.

The Hopf bundle on $ \bS^2$ admits a $U(1)$-action that turns it into an equivariant line bundle, therefore inducing a line bundle $H' \to EU(1) \times_{U(1)} \bS^2$ of fibrewise degree $1$. Let $t:\CPinf \to EU(1) \times_{U(1)} \bS^2$ be the section at $0$. Then the bundle $H:=p^{*} t^{*} (H')^* \otimes H'$ has fibrewise degree $1$ and is trivialised over the zero section.
To finish the proof, let $L:= h_{1}^{*} H$.
\end{proof}

\subsection{Homological stability and its consequences}\label{sec:MarkedPoints}

As we stated in the introduction, proofs of Theorem \ref{modspaceproperties} (\ref{it:stability}) for $n=0$ have appeared in several places \cite{CM,CM2,Boldsen,R-WResolution} with various stability ranges. An explicit reference for the version we are using is \cite[\S 1.4]{R-WResolution}. 

Now we derive the case $n>0$ from the case $n=0$. There is a map
\begin{equation}\label{tangentspacemap}
\epsilon_i:\moddspace{g,r}{n} \lra \CPinf,
\end{equation}
defined as the classifying map of the line bundle $s_{i}^{*} T_v \to \moddspace{g,r}{n}$ (in other words, $\epsilon_i$ takes the tangent space at the $i$-th marked point). The homotopy fibre of $\epsilon_i$ classifies families of surfaces as classified by $\moddspace{g,r}{n}$ but with a framing of $T_v$ at the $i$-th marked point: this is homotopy equivalent to $\moddspace{g,r+1}{n-1}$. Taking all $\epsilon_i$ together, we obtain a fibre sequence
\begin{equation}\label{splitting-and-stripping-sequence}
\modspace{g,r+n} \lra \moddspace{g,r}{n} \stackrel{\epsilon}{\lra} (\CPinf)^n;
\end{equation}
analogous to that used by B{\"o}digheimer and Tillmann in \cite[eq (3.2)]{BT}. Consider the commutative diagram
\[
\xymatrix{
\modspace{g,r+n} \ar[r] \ar[d]& \moddspace{g,r}{n} \ar[r]^-{\epsilon} \ar[d] & (\CPinf)^n \ar@{=}[d]\\
\modspace{h,s+n} \ar[r] &  \moddspace{h,s}{n} \ar[r]^-{\epsilon} &(\CPinf)^n,
}
\]
whose rows are the fibre sequences (\ref{splitting-and-stripping-sequence}) and whose left and middle vertical maps are suitable stabilisation maps. Comparing the Leray--Serre spectral sequences of the two fibre sequences shows homological stability in the case $n>0$. The sequence (\ref{splitting-and-stripping-sequence}) can also be used to prove Theorem \ref{modspaceproperties} (\ref{it:splitting}), using an argument from \cite{BT}. Consider the commutative diagram
\[
\xymatrix{
\modspace{g,r+n} \ar[r] \ar[d]& \moddspace{g,r}{n} \ar[r] \ar[d] & (\CPinf)^n \ar@{=}[d]\\
\modspace{g,r} \ar[r] &  \modspace{g,r} \times (\CPinf)^n \ar[r] &(\CPinf)^n.
}
\]
The bottom sequence is the product sequence, the left vertical map is a composition of stabilisation maps $\gamma(g)$, and the middle vertical map is the product of $\epsilon$ with the forgetful map. Comparing the Leray--Serre spectral sequences of the two fibre sequences establishes that the middle map is a homology isomorphism in degrees $3* \leq 2g-1$.


\subsection{Madsen--Weiss type theorem}\label{sec:MadsenWeiss}

Consider the vector bundle
$$\gamma_{2,n}^\perp \to \mathrm{Gr}^+_2(\bR^{n})$$
given by the orthogonal complement of the tautological bundle, and write $M_n$ for its Thom space. The stabilisation maps $\mathrm{Gr}^+_2(\bR^{n}) \to \mathrm{Gr}^+_2(\bR^{n+1})$ pull back $\gamma_{2, n+1}^\perp$ to $\epsilon^1 \oplus \gamma_{2,n}^\perp$ and so give maps of Thom spaces $\Sigma M_n \to M_{n+1}$. Thus the collection $\{ M_n\}_{n \geq 0}$ forms a \textit{spectrum} in the sense of stable homotopy theory, which is denoted $\MT{SO}{2}$. Similarly, we define a spectrum $\MT{SO}{2} \wedge \CPinf_+$ having as its $n$-th space the smash product $M_n \wedge \CPinf_+$.

Homotopy groups of spectra are very hard to compute in general, but in low degrees, the homotopy groups of $\MT{SO}{2} \wedge \CPinf_+$ may be calculated using the Atiyah--Hirzebruch spectral sequence and the known homotopy groups of $\MT{SO}{2}$ up to degree 3, which have, for example, been computed by the first named author in \cite{EbertIcosahedral}. The result is displayed in Table \ref{table:Homotopy}.
\begin{table}[h]
\caption{}
\centering
\begin{tabular}{lcccccc}
\toprule
$i$ & $-2$ & $-1$ & $0$ & $1$ & $2$ & $3$ \\ \toprule
$\pi_i(\MT{SO}{2})$ & $\bZ$ & 0 & $\bZ$ & 0 & $\bZ$ & $\bZ/24$ \\ 
$\pi_i(\MT{SO}{2} \wedge \CPinf_+)$ & $\bZ$ & 0 & $\bZ^2$ & 0 & $\bZ^3$ & -- \\ \bottomrule
\end{tabular}
\label{table:Homotopy}
\end{table}
The infinite loop space $\Omega^\infty (\MT{SO}{2} \wedge \CPinf_+)$ associated to this spectrum is defined to be the colimit
$$\Omega^\infty (\MT{SO}{2} \wedge \CPinf_+ ):= \colim_{n \to \infty} \Omega^n(M_n \wedge \CPinf_+),$$
and the homotopy groups of this space coincide with those of the spectrum in non-negative degrees. From the table we see that $\Omega^\infty (\MT{SO}{2} \wedge \CPinf_+)$ has $\pi_0$ in bijection with $\bZ^2$. We denote by $\Omega^\infty_0 (\MT{SO}{2} \wedge \CPinf_+)$ the path-component of the basepoint.

Cohen and Madsen \cite{CM, CM2} have defined a comparison map
$$\alpha_{g,k} : \modspace{g,r}(k) \lra \Omega^\infty_{0} (\MT{SO}{2} \wedge \CPinf_+)$$
using Pontrjagin--Thom theory, and shown that it is an integral homology isomorphism in the stable range. Given homology stability, this may also be deduced from the general machines \cite{GMTW, GR-W} for proving Madsen--Weiss type theorems.

Let us now describe how the classes $\mclass_{i,j}$ are defined on the space $\loopinf (\MTSO(2) \wedge \CPinf_+)$. Under the composition
$$\widetilde{H}^{*+2}(BSO(2) \times \CPinf) \cong H^*(\MT{SO}{2} \wedge \CPinf_+) \overset{\sigma_*}\lra H^*(\Omega^\infty_0\MT{SO}{2} \wedge \CPinf_+)$$
of the Thom isomorphism in spectrum cohomology with the cohomology suspension\footnote{For any spectrum $\X = \{X_n\}$, the evaluation maps $\Sigma^n \Omega^n X_n \to X_n$ induce maps on cohomology $H^{*+n}(X_n) \to H^*(\Omega^n X_n)$, which after taking limits over $n$ gives $\sigma^* : H^*(\X) \to H^*(\Omega^\infty \X)$.}, the class $e^{i+1} \otimes c_1^j$ maps to a class we define to be
$$\mclass_{i,j} \in H^{2i+2j}(\Omega^\infty_0\MT{SO}{2} \wedge \CPinf_+;\bZ).$$
The description of push-forwards in terms of Pontrjagin--Thom theory shows that these classes pull back under $\alpha_{g,k}$ to the classes
$$\pi_!(e(T_v)^{i+1} c_1(L)^j) \in H^{2i+2j}(\modspace{g,r}(k);\bZ),$$
as in the introduction.

\subsection{Cohomology of the infinite loop space}\label{sec:CohomologySpectrum}

The rational cohomology of an infinite loop space has a very restricted structure, and is easily deduced from the rational cohomology of the associated spectrum. In turn, the cohomology of a Thom spectrum over a classifying space can be expressed by means of characteristic classes.

\begin{proof}[Proof of Theorem \ref{thm:HolRationalCohomology}]
It is well known that the rational cohomology of the zero component infinite loop space $\Omega_{0}^{\infty} \X$ is given by the free graded-commutative algebra on the vector space $\tau_{* > 0}H^*(\X;\bQ)$ of positive degree elements in spectrum cohomology. (For a proof, note that any rational spectrum splits into a sum of suspensions of $\mathbf{H}\bQ$, and the claim is true for $\Sigma^n \mathbf{H}\bQ$ by direct calculation.)

The rational spectrum cohomology of $\MT{SO}{2}\wedge \CPinf$ in positive degrees is the vector space
$$\bQ \langle u_{-2} \cdot e^{i+1} \wedge c_1^j \,\,|\,\, i+j > 0, j\geq 0, i \geq -1\rangle$$
where $u_{-2} \in H^{-2}(\MT{SO}{2};\bQ)$ denotes the Thom class. The element $u_{-2} \cdot e^{i+1} \wedge c_1^j$ under the cohomology suspension gives the element $\mclass_{i,j}$ on the infinite loop space.
\end{proof}

To prove Theorem \ref{thm:LowDimCohomologyHol}, we first identify $H^2 (\loopinf_0 (\MTSO(2) \wedge \CPinf_+))$ as an abstract group.

\begin{lem}\label{lem:IntegralCohomology}
There are isomorphisms
\[
H^i (\loopinf_0 (\MTSO(2) \wedge \CPinf_+); \bZ)= 
\begin{cases}
0 & i=1,3;\\
\bZ^3 & i=2.
\end{cases}
\]
\end{lem}
\begin{proof}
By Theorem \ref{thm:HolRationalCohomology} we know the statement is true modulo torsion. Thus we must show there is no torsion in cohomology in degrees $* \leq 3$, or equivalently that there is no torsion in homology in degrees $* \leq 2$.

By the isomorphisms listed in Table \ref{table:Homotopy}, $\Omega^\infty_0 (\MT{SO}{2} \wedge \CPinf_+)$ is simply-connected, so in particular $H_1$ is zero and hence torsion free. Furthermore
$$\pi_2(\MT{SO}{2} \wedge \CPinf_+) \cong H_2(\Omega^\infty_0 (\MT{SO}{2} \wedge \CPinf_+);\bZ)$$
by Hurewicz' theorem, and hence we see that the second homology is free abelian of rank three.
\end{proof}

Next we determine a $\bZ$-basis of $H^2(\loopinf_0 (\MTSO(2)\wedge \CPinf_+))\cong \bZ^3$, which begins with naming elements. We describe the elements on the moduli space $\modspace{g}(k)$ as Chern classes of certain index bundles of Cauchy--Riemann operators. A standard exercise in index theory will then show that the classes are indeed induced from $\loopinf(\MTSO(2) \wedge \CPinf_+)$ (this exercise has been solved in the section on ``universal operators" in \cite{EbInd}). We do this in order to deal with more concrete objects and also to facilitate computations.

Recall that $H^2(B \Gamma_g; \bZ)\cong \bZ$ has a generator $\lambda$, the Hodge class, that satisfies $12\lambda = \kappa_1$. To define the next element, we need a divisibility result, based on the Grothendieck--Riemann--Roch theorem. Recall the universal surface bundle $\pi: \mathcal{E}_g(k) \to \modspace{g}(k)$ with vertical tangent bundle $T_v$ and universal line bundle $L \to \mathcal{E}$. Consider the Dolbeault operator $\delbar_{T^{\otimes r} \otimes L^{\otimes s}}$ on the tensor product bundle of $r$ copies of the vertical tangent bundle and $s$ copies of $L$. The Chern character of its index bundle is, by Grothendieck--Riemann--Roch:
$$
\ch (\ind (\delbar_{T^{\otimes r} L^{\otimes s}}))=\pi_{!}(\Td(x) e^{rc_1(T_v)} e^{sc_1(L)})
$$
and the degree $2$ part is
\begin{equation}\label{grothendieckriemannroch}
(6r^2+6r+1)\lambda + \frac{1}{2} s^2(\mclass_{0,1} + \mclass_{-1,2}) + (rs+\frac{1}{2} (s-s^2))\mclass_{0,1}.
\end{equation}
The first and third summand are integral, hence so is the middle summand. In other words, we have proved:
\begin{prop}\label{prop:DivisibilityBy2}
The class $\mclass_{0,1} - \mclass_{-1,2} \in H^2(\modspace{g,r}(k);\bZ)$ is divisible by $2$.
\end{prop}

Thus we may define
$$\zeta := \tfrac{1}{2}(\mclass_{0,1} - \mclass_{-1,2}) \in H^2(\modspace{g,r}(k);\bZ),$$
as this group is torsion-free and hence $\mclass_{0,1} - \mclass_{-1,2}$ is uniquely divisible.

\begin{proof}[Proof of Theorem \ref{thm:LowDimCohomologyHol}]
It is enough to give a map $H^2 (\loopinf_{0} (\MT{SO}{2} \wedge \CPinf_{+}); \bZ) \to \bZ^3$ that maps the tuple $\cB:=(\lambda,\mclass_{0,1},\zeta)$ to a basis. To achieve this, we construct three examples of surface bundles equipped with complex line bundles. The genus and the degree of the line bundle are irrelevant for this purpose. Here are the examples:

\begin{example}\label{example1}
For $g $ in the stable range, we consider the universal surface bundle on $B \Gamma_g$, together with the trivial line bundle on it. We know that $H^2 (B \Gamma_g; \bZ) \cong \bZ \langle \lambda \rangle$, and it is clear that $\cB$ evaluates to $(1,0,0)$.
\end{example}

\begin{example}\label{example2}
Consider the trivial surface bundle $\pi:\CPinf \times \Sigma_g \to \CPinf$, with complex line bundle given by $\pi^* \gamma$. Because this surface bundle is trivial, we have $\kappa_1 =0$, hence $\lambda=0$. Moreover, $\mclass_{0,1}$ is $\chi(\Sigma_g)$ times a generator. Finally, $\mclass_{-1,2}$ is zero. If we put $g=2$, we obtain that $ \cB$ is mapped to $(0,2,1)$.
\end{example}

\begin{example}\label{example3}
Let $H_1$ be the Hirzebruch surface (we use the notation of \cite{MAHirzebruchSurface}). It is an $\bS^2 $-bundle over $\bS^2 $, which is not spin and has $\kappa_1 =0$ (since the signature of the total space is $0$). A basis for $H_2 (H_1; \bZ)$ is given by the fundamental class $u$ of the fibre and the image $v$ of the section $\bS^2 \to H_1$ at $\infty$. Let $(x,y)$ be the Poincar\'e dual basis to $(u,v)$. Using the intersection matrix given in \cite{MAHirzebruchSurface}, it is not hard to see that the Euler class of the vertical tangent bundle is $e=2x+y$. Let $L \to H_1$ be the line bundle with Chern class $y$. Again using the intersection matrix, we compute $\mclass_{0,1}=-1$ and $\mclass_{-1,2}= -1$. Thus $\cB$ is mapped to $(0,-1, 0)$.
\end{example}
\end{proof}

\section{Cohomology of the extended mapping class groups}\label{sec:cohomology-exmcg}

\subsection{Proof of Theorems \ref{thm:ExtendedMCGStability}, \ref{thm:LowDimCohomology} and \ref{thm:PicRationalCohomology}}

The results of the previous section give the stable cohomology of all the extended mapping class groups $\exmcg_{g,r}^{n}$, except for the case $(r,n)=(0,0)$, which is the most interesting of all. Theorem \ref{thm:ExtendedMCGStability} is the key result about the extended mapping class groups in the case $(r,n)=(0,0)$.

\begin{proof}[Proof of Theorem \ref{thm:ExtendedMCGStability}] 
All four spaces in the sequences are connected. By definition the map 
\[\Pi :   \modspace{g}(k)  \lra  \exmcg_{g}(k)\]
induces an isomorphism on fundamental groups. Therefore, by Proposition  \ref{homotopy-type-modspace}, the homotopy fibre of $\Pi$ is a $K(\bZ,2)$.

The mapping space $\map (\CPinf,B \exmcg_g (k))$ is connected (because the source is simply-connected and the target is aspherical). So the map $\Pi \circ \Xi:\CPinf \to B \exmcg_{g}(k)$ is 
nullhomotopic. The choice of a nullhomotopy defines a map $F : \CPinf \to \hofib_*(\Pi)$. To prove that $\Pi$ and $\Xi$ form a fibre sequence, it suffices to prove that $F$ induces an isomorphism on $\pi_2$. But $\hofib_*(\Pi) \to \modspace{g}(k)$ induces an isomorphism on $\pi_2$, so it is enough to prove that $\Xi$ induces an isomorphism on $\pi_2$. However, $\Xi$ is defined as $\CPinf \overset{-\otimes L_0}\to \map(\Sigma_g, \CPinf)(k) \to \modspace{g}(k)$ and both maps are $\pi_2$-isomorphisms. Thus $\Xi$ and $\Pi$ form a fibre sequence.

From now on we replace $\Pi$ by a fibration and we will show that $\exmcg_{g}(k)$ acts trivially on the cohomology of the fibre, using that $g \geq 3$. By Theorem \ref{modspaceproperties} (\ref{it:MW}) and Theorem \ref{thm:LowDimCohomologyHol}, we have the calculation $H_1 (\modspace{g}(k);\bZ)=0$ for $g \geq 3$. The Leray--Serre spectral sequence for the map $\Pi$ implies that $H^1 (B\exmcg_{g}(k);\bZ/2)=0$, and so $\exmcg_g (k)$ can only act trivially on $H^2(\CPinf;\bZ) = \bZ$.

Consequently, the generator of $H^2(\CPinf;\bZ)$ gives an element in the $E_2$-page of the Leray--Serre spectral sequence of the map $\Pi$, which transgresses to a class $\theta \in H^3(B\exmcg_{g}(k);\bZ)$. This corresponds to a homotopy class of map $\theta : B\exmcg_{g}(k) \to K(\bZ, 3)$.
Standard obstruction theory provides a homotopy cartesian square
\begin{equation*}
\xymatrix{
 \modspace{g}(k) \ar[d] \ar[r] & P K(\bZ,3)\simeq \ast \ar[d]\\
  B \exmcg_g(k) \ar[r]^{\theta} & K(\bZ,3)
}
\end{equation*}
(the left-hand vertical map is the path-fibration of $K(\bZ,3)$), which identifies $\modspace{g}(k)$ with the homotopy fibre of $\theta$.

To establish the second part of Theorem \ref{thm:ExtendedMCGStability}, recall that the stabilization map $\modspace{g,1}(k) \to \modspace{g}(k)$ is a homology isomorphism in degrees $3* \leq 2g-1$ by Theorem \ref{modspaceproperties} (\ref{it:stability}). Together with Proposition \ref{homotopy-type-modspace}, this completes the proof.
\end{proof}

We now prepare to prove Theorem \ref{thm:LowDimCohomology}. Recall the fibre sequence
\begin{equation}\label{universalfibration}
\CPinf \overset{\Xi}\lra \modspace{g}(k) \overset{\Pi}\lra B \exmcg_g (k)
\end{equation}
from Theorem \ref{thm:ExtendedMCGStability}. Because $\Pi$ (when replaced by a fibration) is principal, it is simple and so the $E_2$-term of the Leray--Serre spectral sequence takes the form
\[
E^{p,q}_{2} = H^p (B \exmcg_g; H^q (\CPinf)) \cong 
\begin{cases}
H^p (B \exmcg_g;\bZ) & q = 0 \pmod 2\\
0 & q = 1 \pmod 2.
\end{cases} 
\]
The first task will be to compute the edge homomorphism of the spectral sequence, i.e.\ the map $\Xi^* : H^* (\modspace{g}(k)) \to H^* (\CPinf)$. We have described this map geometrically: for $L_0 \to \Sigma_g$ the fixed degree $k$ line bundle, this map classifies the trivial surface bundle $\CPinf \times \Sigma_g \to \CPinf$, together with the line bundle $L=\pr_{1}^{*} \gamma \otimes \pr_{2}^{*} L_0$.

For this surface bundle and line bundle we immediately compute
\begin{align*}
e(T_v) &= (2-2g) \cdot 1 \otimes u\\
c_1(L) &= c_1(\gamma) \otimes 1 + 1 \otimes k\cdot u
\end{align*}
where $u \in H^2(\Sigma_g;\bZ)$ is the Poincar{\'e} dual to a point, from which it is easy to compute the classes $\kappa_{i,j}$.

\begin{lem}\label{lem:EdgeHomomorphism}
The classes $\kappa_{i,j}$ with $i > 0$ map to zero, $\kappa_{0,j}$ maps to $(2-2g) \cdot c_1(\gamma)^j$ and $\kappa_{-1,j}$ maps to $j \cdot k \cdot c_1(\gamma)^{j-1}$.

In particular, $\lambda$ goes to zero, $\zeta$ goes to $(1-g-k) \cdot c_1(\gamma)$ and $\mclass_{0,1}$ to $(2-2g) \cdot c_1(\gamma)$.
\end{lem}

\begin{proof}[Proof of Theorem \ref{thm:LowDimCohomology}]
We study the Leray--Serre spectral sequence in cohomology for the fibre sequence (\ref{universalfibration}), which has $E_2$ page as shown in Figure \ref{fig:SerreSS}
\begin{figure}[h]
\centering
\includegraphics[bb=0 0 195 126]{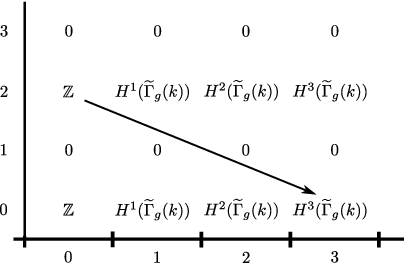}
\caption{$E_2$ page of the Serre spectral sequence for (\ref{universalfibration}).}
\label{fig:SerreSS}
\end{figure}
\noindent and which converges to zero in total degrees 1 and 3 by Lemma \ref{lem:IntegralCohomology}. 
From this we deduce that $H^1(B \exmcg_{g}(k))=0$ and we get an exact sequence
\[
0 \lra H^2 (B \exmcg_{g}(k)) \lra H^2 (\modspace{g}(k)) \lra H^2 (\CPinf) \stackrel{d_3}{\lra} H^3 (B \exmcg_{g}(k)) \lra 0.
\]
By Lemma \ref{lem:EdgeHomomorphism}, the index of the image of the edge homomorphism in $H^2 (\CPinf)$ is $\gcd(2-2g, 1-g-k)$. The kernel of the edge homomorphism is free of rank $2$, and by elementary algebra, it is generated by $\lambda$ and $\eta$.
\end{proof}

We nor prepare to prove Theorem \ref{thm:PicRationalCohomology}, where we describe $H^* (B \exmcg_{g}(k); \bQ)$ in the stable range. We must begin with the definition of the classes $\nclass_{i,j}$. Let $B\Gamma^{1}_{g} \to B\Gamma_g$ be the universal surface bundle. There is a homotopy commutative diagram
\begin{equation*}
\xymatrix{
\CPinf \ar@{=}[d] \ar[r]^{\iota} & \mathcal{E}_g(k) \ar[d]^{\pi} \ar[r]^-{f}  &  B \exmcg_{g}^{1}(k)  \ar[d]^{p} \ar[r] & B \Gamma_{g}^{1} \ar[d]\\
\CPinf \ar[r]^{\epsilon}  & \modspace{g}(k) \ar[r]^-{\Pi} & B \exmcg_g (k) \ar[r] & B\Gamma_g.
}
\end{equation*}

The two squares on the right are cartesian and $\pi$ and $p$ are surface bundles. The maps $\iota$ and $\epsilon$ are the inclusions of the respective homotopy fibres.
The first two maps of the bottom sequence form a fibre sequence by Theorem \ref{thm:ExtendedMCGStability}.
Recall that $\mathcal{E}_g = \map(\Sigma_g^1,\CPinf) \times_{\Difff{g}{1}} E\Difff{g}{1}$. 
Evaluation at the basepoint of $\Sigma_g^1$ is a $\Difff{g}{1}$-equivariant map $\map(\Sigma_g^1,\CPinf) \to \CPinf$ 
and so defines a map $\mathcal{E}_{g} \to \CPinf$, which is a left homotopy inverse to $\iota$.


As $\iota$ has a left homotopy inverse, there is a short exact sequence
$$0 \lra H^2 (\exmcg_{g}^{1}(k);\bQ) \lra H^2 (\mathcal{E}_g(k);\bQ) \lra H^2 (\CPinf;\bQ) \lra 0.$$
The rational cohomology class $c_1(L) - \tfrac{1}{2-2g}\pi^* \mclass_{0,1}$ maps to $0 \in H^2 (\CPinf;\bQ)$ by Lemma \ref{lem:EdgeHomomorphism}, and thus there exists a unique $v \in H^2 ( \moddspace{g}{1}(k);\bQ)$ with $f^* v = c_1(L) - \tfrac{1}{2-2g}\pi^* \mclass_{0,1}$. Let $e(T_v) \in H^2 (\moddspace{g}{1}(k);\bZ)$ be the Euler class of the vertical tangent bundle of $p$, and define
\begin{equation*}
\nclass_{i,j} := p_{!}(e(T_v)^{i+1}v^j) \in H^{2i+2j}(B \exmcg_{g}(k);\bQ).
\end{equation*}
We compute
\begin{align*}
\Pi^* \nclass_{i,j} &= \pi_{!} (e(T_v)^{i+1} (c_1(L) - \tfrac{1}{\chi}\pi^* \mclass_{0,1})^j)\\
&= \sum_{p+q=j} \binom{j}{p} (-\tfrac{1}{\chi})^q \pi_{!}(e(T_v)^{i+1} c_1(L)^p \pi^{*} (\mclass_{0,1}^q))\\
&= \sum_{p+q=j} \binom{j}{p} (-\tfrac{1}{\chi})^q \mclass_{i,p}\mclass_{0,1}^q,
\end{align*}
bearing in mind that the degree zero classes are $\mclass_{-1, 1} = k$, $\mclass_{0,0} = 2-2g$.

\begin{proof}[Proof of Theorem \ref{thm:PicRationalCohomology}]
To see that $\Pi^*$ is rationally injective, we use the fibre sequence (\ref{universalfibration})
and the fact from Theorem \ref{thm:LowDimCohomology} that the element $\theta \in H^3(B\exmcg_g(k);\bZ)$ which classifies this fibre sequence is torsion. This implies that the rational Leray--Serre spectral sequence collapses, so that $\Pi^*$ is rationally injective.

We certainly have an algebra homomorphism
\begin{equation}\label{eq:freealg}
\bQ[\nclass_{i,j} \mid i\geq -1, j \geq 0, i+j >0; (i,j)\neq (0,1)] \lra H^{*}(B \exmcg_{g}(k); \bQ)
\end{equation}
and by the calculation above $\Pi^* \nu_{i,j} \equiv \kappa_{i,j} \; \mathrm{mod} \;(\kappa_{0,1})$, so the algebra $\bQ[\nclass_{i,j} \mid i\geq -1, j \geq 0, i+j >0; (i,j)\neq (0,1)]$ maps injectively into $H^*(\modspace{g}(k);\bQ)$ in the stable range, and so (\ref{eq:freealg}) is injective too. Counting dimensions and using the collapse of the Leray--Serre spectral sequence above shows (\ref{eq:freealg}) is also surjective in the stable range.
\end{proof} 

\begin{rem}
The computation above shows that $\nclass_{0,1}=0$ and $\nclass_{-1,2}= \tfrac{1}{g-1}\gcd (2g-2,g-1+k) \eta$.
\end{rem}

\subsection{Structure of the extended mapping class groups}

From the results proved so far and some facts on the cohomology of the mapping class groups, we can illuminate the structure of the extended mapping class groups. If we write $H_{g,r}^{n}:=H^1 (\Sigma_{g,r}^{n}, \partial \Sigma_{g,r}^n \cup P;\bZ)$, then there is a homotopy equivalence $\map_\partial (\Sigma_{g,r}^n, \CPinf)(k) \simeq BH_{g,r}^n$ as long as $r+n > 0$, and so extensions
\begin{equation}\label{groupextension}
0 \lra H_{g,r}^{n} \lra \exmcg_{g,r}^{n}(k) \lra \Gamma_{g,r}^{n} \lra 1.
\end{equation}
If $(r,n)=(0,0)$ there is a square of isomorphisms
\begin{equation}\label{eq:action}
\begin{gathered}
\xymatrix{
\pi_1(\map(\Sigma_g, \CPinf)(k), L_0) \ar[d]^{\text{Hurewicz}} & \ar[d]^{\text{Hurewicz}} \ar[l]_{-\otimes L_0}  \pi_1(\map(\Sigma_g, \CPinf)(0), *) \cong H_g\\
H_1(\map(\Sigma_g, \CPinf)(k);\bZ) & \ar[l]_{-\otimes L_0}  H_1(\map(\Sigma_g, \CPinf)(0);\bZ)
}
\end{gathered}
\end{equation}
and so we still have the sequence (\ref{groupextension}).

\begin{lem}
The action of $\Gamma_{g,r}^n$ on $H_{g,r}^n$ in the extension (\ref{groupextension}) is the usual one.
\end{lem}
\begin{proof}
Let $\gamma \in \Difff{g,r}{n}$ be a diffeomorphism. It acts on $\map_\partial(\Sigma_{g,r}^n, \CPinf)(k)$ by precomposition, which does not typically fix the point $L_0$. Instead, we have the equation $\gamma^* \circ (- \otimes L_0) = (-\otimes \gamma^*L_0) \circ \gamma^*$. As $\gamma$ has degree 1, $\gamma^* L_0$ is homotopic to $L_0$, and so the action of $\gamma$ on $H_1(\map(\Sigma_g, \CPinf)(k);\bZ)$ corresponds to the usual action on $H_g$ under the isomorphism of (\ref{eq:action}).
\end{proof}

We ask when the extension (\ref{groupextension}) is split, i.e.\ admits a homomorphic section.
For all values of $r,n,k$ in which Proposition \ref{existencelinebundle} applies, there is a splitting, showing triviality in these cases. For $(r,n)=(0,0)$, the following proposition shows that this result is sharp (in the stable range).

\begin{prop}
Let $g \geq 6$ and $n=r=0$. Then the extension (\ref{groupextension}) splits if and only if $k \equiv 0 \pmod {2g-2}$.
\end{prop}

\begin{proof}
Let $\sigma$ be a splitting, giving a map $\sigma:B \Gamma_{g} \to B \exmcg_g (k)$. If $g \geq 6$, then it is known that $H^3 (B \Gamma_g;\bZ)=0$ and so $\sigma^* \theta =0$. Therefore $\sigma$ lifts to a cross-section $\tau: B \Gamma_g \to \modspace{g}(k)$ by elementary obstruction theory. This means that there is a line bundle on the universal surface bundle $E_{g} \to B \Gamma_g$ of fibrewise degree $k$. But B\"odigheimer and Tillmann have shown \cite{BT} that the fibrewise degree of a line bundle on $E_g$ is divisible by the Euler number $2-2g$.
\end{proof}

\begin{rem}
Morita \cite{Morita86} has proved that $H^2 (\Gamma_g; H_g) \cong \bZ/(2-2g)$. One can show that the map $\bZ \to H^2 (\Gamma_g; H_g)$ that sends $k$ to the isomorphism class of the extension $H_g \to \exmcg_g (k) \to \Gamma_g$ is a surjective homomorphism.
\end{rem}

\subsection{Proof of Theorem \ref{thm:DegenerationSpectralSequence}}

The rational Leray--Serre spectral sequence of the fibre sequence
\begin{equation}\label{homology-exmcg-mcg-fibration}
B H_{g,r}^{n} \lra B \exmcg_{g,r}^{n} (k) \lra B \Gamma_{g,r}^{n}
\end{equation}
has the form
\begin{equation}\label{specseq}
E_{2}^{p,q} =  H^p (B \Gamma_{g,r}^{n}; \wedge^q (H_{g,r}^{n} \otimes \bQ)) \Longrightarrow H^{p+q} (B \exmcg_{g,r}^{n}; \bQ),
\end{equation}
where the action of $\Gamma_{g,r}^{n}$ is given by the usual action on the cohomology of $\Sigma_{g,r}^{n}$. The proof that this collapses begins with the case $(r,n)=(0,0)$, where we need the following result.

\begin{prop}\label{thm:deligne}(Deligne \cite[Proposition 2.1]{Deligne})
Let $f:X \to Y$ be a fibration, $\cF$ be a local coefficient system of $\bQ$-vector spaces on $X$, $m \in \bN$ and $u \in H^2 (X; \bQ)$ be given. Denote $X_y:= f^ {-1}(y)$ and suppose that for each $y \in Y$ and each $i \geq 0$, the multiplication by $(u^i)|_{X_y}$ induces an isomorphism $H^{m-i} (X_y; \cF|_{X_y}) \to H^{m+i} (X_y; \cF|_{X_y})$ (the ``Lefschetz condition for $(X,f,\cF,u,m)$"). Then the Leray--Serre spectral sequence $E_{2}^{p,q}\Rightarrow H^{p+q}(X;\cF)$ for $f$ collapses at the $E_2$-stage.
\end{prop}

To apply this criterion, let $\omega \in H^2 (BH_g)$ be the symplectic class. It is easy to see that multiplication by $\omega^k$ induces an isomorphism $H^{g-k}(BH_g;\bQ) \to H^{g+k}(BH_g; \bQ)$. The next lemma therefore shows that (\ref{homology-exmcg-mcg-fibration}) satisfies the Lefschetz condition for the constant coefficient system $\bQ$, $m=g$ and the class $\eta \in H^2 (B \exmcg_g (k))$ from Theorem \ref{thm:LowDimCohomology}. 

\begin{lem}\label{kaehlerclass}
The restriction of $\eta \in H^2 (B \exmcg_g (k))$ to $B H_g$ is a non-trivial multiple of $\omega$.
\end{lem}

\begin{proof}
As the $\Gamma_g$-invariant part of $H^2(BH_g;\bZ) = \wedge^2 H_g$ is the subspace spanned by the symplectic form $\omega$, the restriction of $\eta$ is certainly some multiple of $\omega$. By a theorem of Morita \cite[Proposition 4.1]{Morita89}, $H^1(\Gamma_g;H_g)=0$. Thus from the spectral sequence (\ref{specseq}) we find an exact sequence
$$0 \lra H^2(\Gamma_g;\bQ) \lra H^2(\exmcg_g(k);\bQ) \lra H^2(H_g;\bQ)^{\Gamma_g} \overset{d_3}\lra \cdots.$$
As $\eta$ and $\kappa_{1,0}$ are linearly-independent, the image of $\eta$ is non-trivial.
\end{proof}

\begin{rem}
In fact, one may calculate that the restriction of $\eta$ is $\frac{2-2g}{\gcd(2-2g, g+k-1)}\cdot \omega$.
\end{rem}

This shows that the spectral sequence (\ref{specseq}) collapses if $r=n=0$.

\begin{proof}[Proof of the collapse of (\ref{specseq}) for $r+n=1$ or $r>1$ and $n=0$]
If $g:Z \to X$ is a map, then the Lefschetz condition for a fibration $X \to Y$ implies that for the pullback fibration $g^* X \to Z$. Because the two squares
\[
\xymatrix{
B \exmcg_{g,1}(k) \ar[r] \ar[d]&B \exmcg_{g}^{1}(k) \ar[r] \ar[d]&  B \exmcg_g(k) \ar[d] \\
B \Gamma_{g,1} \ar[r] &B \Gamma_{g}^{1} \ar[r] &  B \Gamma_g
}
\]
are cartesian, this observation shows the collapse of the spectral sequence for $r+n=1$. 

Now we use homological stability to extend this result to the cases $n=0$ and $r\geq 1$, so we only obtain the collapsing in the stable range. We use the homological stability theorem of Boldsen \cite{Boldsen} for the mapping class group with twisted coefficients. The coefficient system $\wedge^j H_\bQ$ has degree $j$ and so the map
$$H^{2i+j}(\Gamma_{g,r};\wedge^j H_\bQ) \lra H^{2i+j}(\Gamma_{g,1};\wedge^j H_\bQ)$$
is an isomorphism for $2i+j \leq \lfloor\tfrac{2g}{3}\rfloor-j$. Hence the map of spectral sequences corresponding to $\exmcg_{g,1}(k) \to \exmcg_{g,r}(k)$ is an isomorphism in the stable range on $E_2$, and hence the spectral sequence for $\exmcg_{g,r}(k)$ also collapses.
\end{proof}

\begin{proof}[Computation of the associated graded algebras for $n=0$, $r>0$]
The case $r=1$, $k=0$ has been proved by Kawazumi \cite{Kawazumi98, KawazumiInfinitesimal, KawazumiGoD}. He defined certain cohomology classes \cite{Kawazumi98} 
$$\widetilde{m}_{i,j} \in H^{2i+2j-2}(\exmcg_{g,1}(0);\bZ),$$
in the framework of group cohomology, as pushforwards $\pi_!(e^i \cdot \widetilde{\omega}^j)$ where the element $e \in H^2({\Gamma}_{g,1}^1, {\Gamma}_{g,1} \times \bZ;\bZ)$ is the group level analogue of the Euler class of the vertical tangent bundle, and $\widetilde{\omega} \in H^2(\exmcg_{g,1}^1, \exmcg_{g,1} \times \bZ;\bZ)$ is represented by a cocycle he defines manually.  He does not phrase it in this way, but it is the Euler class of the relative central extension
\begin{equation}\label{eq:KawazumiCentralExtension}
\begin{gathered}
\xymatrix{
\bZ \ar[r] \ar@{=}[d] & H_1(\Sigma_{g,1}^1) \rtimes \Gamma_{g,1}^1  \ar[r] & H_1(\Sigma_{g,1}) \rtimes \Gamma_{g,1}^1 \\
\bZ \ar[r] & (H_1(\Sigma_{g,1}^1) \rtimes \Gamma_{g,1}) \times \bZ \ar[r] \ar[u] & \exmcg_{g,1}(0) \times \bZ \ar[u].}
\end{gathered}
\end{equation}
The homomorphism $\rho:\Gamma_{g,1} \times \bZ \to \Gamma_{g,1}^{1}$ implicit in the diagram is obtained by gluing in an annulus with a Dehn twist, see \cite[p.\ 140]{Kawazumi98}.
Here $H_1(\Sigma_{g,1}^1)$ denotes the homology of the surface with the marked points removed. The top sequence arises from the short exact sequence $0 \to \bZ \to H_1(\Sigma_{g,1}^1) \to H_1(\Sigma_{g,1}) \to 0$ of $\Gamma_{g,1}^{1}$-modules. The bottom sequence is induced from the top one via $\rho$ and it is split via the $\Gamma_{g,1}$-equivariant inclusion $H_1(\Sigma_{g,1}) \to H_1(\Sigma_{g,1}^1)$ obtained by boundary connected sum with a punctured disc.

Diagram (\ref{eq:KawazumiCentralExtension}) has a topological interpretation. Recall the $U(1)$-principal bundle $\moddspace{g,1}{1} \to \mathcal{E}_{g,1}$ (that is, the frame bundle of the universal line bundle) and recall also that over the boundary $\modspace{g,1} \times \bS^1$ this bundle is trivialised. The diagram
\begin{equation*}
\xymatrix{
\bS^1 \ar[r] \ar@{=}[d]& \moddspace{g,1}{1} \ar[r] & \mathcal{E}_{g,1} \\
\bS^1 \ar[r] & \partial \mathcal{E}_{g,1} \times \bS^1 \ar[u] \ar[r] & \partial \mathcal{E}_{g,1} \ar[u] \ar@{=}[r]& B \exmcg_{g,1}(0) \times \bS^1 
}
\end{equation*}
consisting of aspherical spaces induces the diagram (\ref{eq:KawazumiCentralExtension}) on fundamental groups. Therefore $\widetilde{\omega}$ coincides with the Chern class $c_1 (L) $. Thus we have proved:

\begin{prop}
There is an equality $\widetilde{m}_{i,j} = \mclass_{i-1,j} \in H^*(\exmcg_{g,1}(0);\bZ)$.
\end{prop}

We say that the \emph{natural filtration} of the cohomology of $\exmcg_{g,r}^n$ is the Leray filtration coming from the fibre sequence (\ref{homology-exmcg-mcg-fibration}).
Kawazumi shows that the class $\widetilde{m}_{i+1,j} \in H^* (B\exmcg_{g,1})$ is detected in the Leray--Serre spectral sequence (\ref{specseq}) by a certain class
$$m_{i+1, j} \in H^{2i+j}(\Gamma_{g,1};\wedge^j H)$$
which he constructs. Equivalently, we may say that $\mclass_{i,j}=\widetilde{m}_{i+1,j}$ has natural filtration precisely $j$. In \cite[Theorem 3.2]{KawazumiGoD} he then shows that rationally the associated graded algebra to the natural filtration is
$$\mathcal{G}r(H^*(\exmcg_{g,1}(0);\bQ)) = H^*(\Gamma_{g,1};\bQ) \otimes \bQ[m_{i+1,j} \, |\, j > 0],$$
which proves Theorem \ref{thm:DegenerationSpectralSequence} in this case. On the other hand, by Theorem \ref{modspaceproperties} (\ref{it:independence}), there is an equivalence $B\exmcg_{g,1}(k) \simeq B\exmcg_{g,1}(0)$ so this proves the result for all $k$. Furthermore, homological stability for twisted coefficients extends this to all $r > 0$ by the comparison maps $B\exmcg_{g,r}(k) \to B\exmcg_{g,1}(k)$.
\end{proof}

\begin{proof}[Proof of collapsing and associated graded algebra for $n>0$ and $r \geq 0$.]
As in Section \ref{sec:MarkedPoints} we have a fibre sequence
$$B\Gamma_{g, r+n} \lra B\Gamma_{g,r}^n \overset{e_1 \times \cdots \times e_n}\lra \CPinf \times \cdots \times \CPinf$$
and as long as $r > 0$ the composition $B\Gamma_{g, r+n} \to B\Gamma_{g,r}^n \to B\Gamma_{g, r}$ induces a cohomology isomorphism with coefficients in $\wedge^j H_\bQ$ in a stable range. Thus
$$H^*(\Gamma_{g,r};\wedge^j H_{\bQ})[e_1, ..., e_n] \lra H^*(\Gamma_{g,r}^n;\wedge^j H_{\bQ})$$
is an isomorphism in a stable range. The same argument with the fibre sequence 
$$B\exmcg_{g, r+n}(k) \lra B\exmcg_{g,r}^n(k) \overset{e_1 \times \cdots \times e_n}\lra \CPinf \times \cdots \times \CPinf$$
shows that $H^*(\exmcg_{g,r}(k);\bQ)[e_1, ..., e_n] \to H^*(\exmcg_{g,r}^n(k);\bQ)$ is also an isomorphism in the stable range, and by counting dimensions we see that the spectral sequence for $\exmcg_{g,r}^n(k)$ must collapse.

This leaves the case $r=0$, $n > 0$. We consider the fibre sequence
$$B{\Gamma}_{g, 1}^{n-1} \lra B{\Gamma}_{g}^n \overset{e_1}\lra \CPinf$$
which gives a spectral sequence
$$\hat{E}_2^{*,*} := H^*(\Gamma_{g,1}^{n-1}; \wedge^q H_{\bQ})[e_1] \Longrightarrow H^*(\Gamma_g^n;\wedge^qH_\bQ).$$
The output of this spectral sequence is the input of the spectral sequence
$$E_2^{p,q} := H^p(\Gamma_g^n;\wedge^q H_{\bQ}) \Longrightarrow H^{p+q}(\exmcg_g^n(k);\bQ)$$
which we want to show collapses. In the stable range
$$H^*(\exmcg_g^n(k);\bQ) \cong \bQ[e_1, ..., e_n, \kappa_1, \kappa_2, ...] \otimes \bQ[\mclass_{i,j} \,\, | \,\, i+j>0, j>0, i \geq -1]$$
and
$$\bigoplus_{p, q} H^p(\Gamma_{g,1}^{n-1}; \wedge^q H_{\bQ}) \cong \bQ[e_2, ..., e_n, \kappa_1, \kappa_2, ...] \otimes \bQ[x_{i,j} \,\, | \,\, i+j>0, j>0, i \geq -1]$$
and hence by counting dimensions we see that both spectral sequences must collapse. From the collapse of the first spectral sequence we find that the associated graded is as claimed.
\end{proof}

\begin{proof}[The case $(r,n)=(0,0)$]
The last part of Theorem \ref{thm:DegenerationSpectralSequence} to be proved is the claim about the associated graded to the natural filtration of $H^*(\exmcg_g (k);\bQ)$.

For this we consider the fibre sequences
$$\Sigma_g \lra B\Gamma_g^1 \overset{p}\lra B\Gamma_g \quad\quad \Sigma_g \lra B\exmcg_g^1 (k) \overset{\tilde{p}}\lra B\exmcg_g (k).$$
The first one satisfies the assumptions of Theorem \ref{thm:deligne} for any rational coefficient system that is induced from $B \Gamma_g$, by taking $u$ to be the Euler class of the vertical tangent bundle. The second sequence is the pullback of the first.
Therefore the spectral sequence for $B \Gamma_{g}^{1}\to B \Gamma_g$ with coefficients in $\wedge^r H_{\bQ}$:
\[
E_{2}^{p,q} = H^p (B \Gamma_g; H^q (\Sigma_g) \otimes \wedge^r H_{\bQ}) \Longrightarrow H^{p+q}(B \Gamma_{g}^{1}; \wedge^r H_{\bQ})
\]
collapses at the $E_2$-stage and 
$$p^* : H^*(\Gamma_g;\wedge^q H_{\bQ}) \lra H^*(\Gamma_g^1;\wedge^q H_{\bQ})$$
is injective. Hence the natural filtration on $H^*(\exmcg_g(k);\bQ)$ and the filtration induced by the injection $\tilde{p}^*$ agree, and the claimed description of the associated graded follows from Theorem \ref{thm:PicRationalCohomology}.
\end{proof}

\section{Relation to complex algebraic geometry}\label{sec:Geometry}

\subsection{Definitions and results}

We first summarise our results, and will give proofs below. We assume that the reader is familiar with the basic vocabulary of the language of topological stacks. If not, they are advised to read \cite{Hein}, \cite{Noo} and the relevant parts of \cite{EbGian}.

\begin{defn}
Let $\HolStack_g$ denote the following stack, defined on the site $\Top$. An object of $\HolStack_g(X)$, $X\in \Top$, is a triple $(E, \pi,L)$, consisting of a family of genus $g$ Riemann surfaces $\pi : E \to X$ and a fibrewise holomorphic line bundle $L \to E$. An isomorphism $(E, L) \to (E', L')$ is an isomorphism of families $f : E \to E'$ and an isomorphism $h : L \to f^*L'$ of holomorphic line bundles. We denote by $\HolStack_{g}^{k} \subset \HolStack_g$ the open and closed substack consisting of those bundles having fibrewise degree $k$.
\end{defn}

There is an obvious forgetful map $\HolStack_g \to \ModStack_g$, which just remembers the underlying family of Riemann surfaces.

The definition of the Picard stack $\PicStack_g$ is a little more involved. Let $\pi : E \to X$ be a family of genus $g$ Riemann surfaces, that is, an element of $\ModStack_g (X)$. We define the associated Picard bundle $\Pic(E/X) \to X$ as follows. Let $\mathcal{O}$ denote the sheaf of continuous, fibrewise holomorphic functions on $E$, and $\mathcal{O}^\times$ the subsheaf of nowhere zero functions. The exponential sequence of sheaves $\bZ \to \mathcal{O} \to \mathcal{O}^\times$ gives an exact sequence of sheaves on $X$,
\begin{equation*}
0 \lra R^1 \pi_*\bZ \lra R^1\pi_*\mathcal{O} \lra R^1\pi_*\mathcal{O}^\times \lra R^2 \pi_*\bZ \lra 0.
\end{equation*}

This is because the sequence $0 \to \pi_* \bZ \to \pi_* \cO \to \pi_* \cO^{\times} \to 0$ is isomorphic to $0 \to \bZ \to \bC \to \bC^{\times} \to 0$ and hence exact. Moreover, $R^2 \pi_* \cO=0$.
These are all sheaves of continuous sections of certain bundles of groups on $X$,
\begin{equation*}
0 \lra [H^1(F_x;\bZ)] \lra [H^1(F_x;\mathcal{O})] \lra \Pic(E/X) \lra \bZ \times X \lra 0
\end{equation*}
where $[H^1(F_x;\mathcal{F})]$ is the bundle of fibrewise first cohomologies with coefficients in the sheaf $\mathcal{F}$, and $\Pic(E/X)$ is a bundle of abelian groups isomorphic to $\bZ \times \bT^{2g}$. The fibre of $\Pic(E/X)$ over $x \in X$ is the Picard group $\Pic(\pi^{-1}(x))$. The group $\pi_1(X)$ acts trivially on the set of path components, and we denote by $\Pic^k(E/X)$ the degree $k$ component.

\begin{defn}
Let $\PicStack_g$ denote the following stack, defined on the site $\Top$. An object of $\PicStack_g(X)$ is a family of Riemann surfaces $\pi : E \to X$ and a section $s : B \to \Pic(E/X)$. An isomorphism $(E, s) \to (E',s')$ is an isomorphism of families $f : E \to E'$ such that $f^*(s') = s$. We denote by $\PicStack_g^k \subset \PicStack_g$ the open and closed substack consisting of those pairs $(E, s)$ where $s$ takes values in $\Pic^k(E) \subset \Pic(E)$.
\end{defn}

\begin{defn}
Let $\Phi_g^k : \HolStack_g^k \to \PicStack_g^k$ be the morphism that sends a fibrewise holomorphic line bundle to its isomorphism class.
\end{defn}

Let $X$ be a space and $\beta:X \to \PicStack_g^k$ be a map, which is a family of Riemann surfaces $\pi:E \to X$ and a section $s: X \to \Pic^k (E/X)$. We can ask whether there exists a (fibrewise) holomorphic line bundle $L \to E$ such that for each $x \in X$, the isomorphism class of the line bundle $L|_{E_x}$ is $s(x)$. Such a line bundle will be called \emph{topological Poincar\'e line bundle} for $\beta$. It is plain to see that topological Poincar\'e line bundles for $\beta$ are the same as lifts $X \to \HolStack_g^k$ of $\beta$ along $\Phi_g^k$. This observation can be used to generalise the notion of a Poincar\'e line bundle to stacks, as follows: a Poincar\'e line bundle for a map $\mathsf{X} \to \PicStack_{g}^ {k}$ is a lift to $\HolStack_{g}^{k}$.

\begin{rem}
In the classical literature, the notion ``Poincar\'e line bundle" is used slightly differently. Namely, for a single Riemann surface $S$, one asks for a line bundle $L \to \Pic^k (S) \times S$, such that for each $\ell \in \Pic^k (S)$, the restriction $L|_{\{\ell\} \times S}$ is in the isomorphism class $\ell$. More generally, for a family $E \to X$, one asks for a \emph{holomorphic} line bundle $L \to E \times_{X} \Pic^k (E/X)$ with that property.
\end{rem}

\begin{thm}\label{thm:Gerbe}
Let $\beta: X \to \PicStack_g^k$ be given. Then the following properties hold:
\begin{enumerate}[(i)]
\item Each $x \in X$ has an open neighborhood $U$ such that the restriction $\beta|_U$ admits a Poincar\'e line bundle.
\item If $L_i$, $i=0,1$ are two Poincar\'e line bundles, then each $x \in X $ has a neighborhood $U$, such that $(L_0)|_U \cong (L_1)|_U$.
\item If $L$ is a Poincar\'e line bundle, then the group of automorphisms of $L$ is isomorphic to $\map(X, \bC^{\times})$.
\end{enumerate}
More technically, the map $\Phi_g^k$ is a gerbe with band $\bC^\times$.
\end{thm}

In fact, the three statements in Theorem \ref{thm:Gerbe} are precisely the gerbe axioms \cite[\S V.2 ]{Bryl}, translated to the present context. Probably Theorem \ref{thm:Gerbe} is a special case of a statement that is well-known among algebraic geometers who are also fond of stacks, though we were unable to find a complete proof in the literature. The proof below applies only to line bundles on curves of genus $g \geq 2$. 

A map of stacks $\mathsf{X} \to \mathsf{Y}$ of stacks is a \emph{universal weak equivalence} if for each space $Z$ and each map $Z \to \mathsf{Y}$, the induced map $\mathsf{X} \times_{\mathsf{Y}} Z \to Z$ is a weak equivalence of topological spaces. A \emph{homotopy type} of the stack $\mathsf{X}$ is a universal weak equivalence $X \to \mathsf{X}$ with source a space $X$. Each topological stack admits a homotopy type, and two homotopy types are unique up to homotopy equivalence. See Noohi \cite{Noo} for these notions (he writes ``classifying space" instead of ``homotopy type"). We say shortly ``the homotopy type of the stack $\mathsf{X}$ is $X$" if there exists a space $X'$, a weak equivalence $X' \simeq X$ and a universal weak equivalence $X' \to \mathsf{X}$.

\begin{thm}\label{thm:HolCxStacks}
The homotopy type of the stack $\HolStack_g^k$ is $\modspace{g}(k)$, the homotopy type of $\PicStack_g^k$ is $B \exmcg_g(k)$, and under these equivalences the gerbe $\Phi_g^k$ corresponds to $\Pi$.
\end{thm}

If $\mathsf{Y} \to \mathsf{X}$ is a gerbe with band $\bC^{\times}$, the map $X \to Y$ on homotopy types has homotopy fibre a $K(\bZ,2)$ at every point, and the associated fibration is principal (c.f.\ \cite{Noo2}).


Isomorphism classes of line bundles on a stack are in bijection with the sheaf cohomology group $H^1 (\mathsf{X}; \underline{\bC^{\times}})$. Similarly, isomorphism classes of gerbes with band $\bC^{\times}$ are in bijection with the sheaf cohomology group $H^2 (\mathsf{X}; \underline{\bC^{\times}})$. The connecting homomorphism $H^i (\mathsf{X}; \underline{\mathbb{C}^{\times}}) \to H^{i+1}(\mathsf{X}; \bZ)$ may \emph{fail} to be an isomorphism on a general stack. However, if $\mathsf{X}$ is a \emph{local quotient stack} (see \cite{Hein} or \cite{EbGian} for this notion), the connecting homomorphism is an isomorphism. The reason is that the sheaf of continuous functions on a local quotient stack is acyclic. This is well-known, but we do not know a reference and include a proof as Lemma \ref{lem:TopPicardGpInj}.

The image of the isomorphism class of a gerbe $\cG$ on $\mathsf{X}$ in $H^3 (\mathsf{X};\bZ)$ is called \emph{Dixmier--Douady class}. So on a local quotient stack, gerbes are classified by their Dixmier--Douady class. Finally, the Dixmier--Douady class of $\cG$ agrees with the characteristic class of the induced principal fibration. 

\begin{prop}\label{local-quotient-stacks}
Both $\HolStack_g$ and $\PicStack_g$ are local quotient stacks.
\end{prop}

Together with the observation that a gerbe is trivial once there exists a section, we obtain a corollary, using Theorem \ref{thm:LowDimCohomology}.

\begin{cor}
There exists a Poincar\'e line bundle for the identity map on $\PicStack_g^k$ if and only if $\gcd(2g-2,1-g-k)=1$.
\end{cor}

This is a precise topological analogue with the result by Mestrano and Ramanan \cite{MestranoRamanan} that takes place in the holomorphic category.
The proof of Proposition \ref{local-quotient-stacks} will show the stronger result that both stacks are even holomorphic local quotient stacks. On a holomorphic stack $\mathsf{X}$, we can talk about holomorphic line bundles. Let $\Pic_{\hol}(\mathsf{X})$ be the group of isomorphism classes of holomorphic line bundles. There is a comparison map $\Pic_{\hol} (\mathsf{X})\to \Pic_{\top} (\mathsf{X})$, which we will show to be an isomorphism in Theorem \ref{comparison:linebundles}. This theorem, combined with the purely topological Theorem \ref{thm:LowDimCohomology} and Lemma \ref{kaehlerclass}, gives another proof of the main result of Kouvidakis' paper \cite{Kouvidakis}.

\subsection{The gerbe of holomorphic maps}\label{sec:gerbe}

Here we prove Theorem \ref{thm:Gerbe}, which begins with a trivial observation:

\begin{lem}\label{periodicity}
Tensor multiplication with the cotangent bundle induces a commutative diagram
\begin{equation*}
\xymatrix{
\HolStack_g^k \ar[r] \ar[d]^{\Phi_{g}^{k}} & \HolStack_{g}^{k+2g-2} \ar[d]^{\Phi_{g}^{k+2g-2}} \\
\PicStack_g^k \ar[r] & \PicStack_{g}^{k+2g-2}
}
\end{equation*}
whose horizontal arrows are isomorphisms.
\end{lem} 

Therefore, it is enough to prove Theorem \ref{thm:Gerbe} for large values of $k$. In the sequel, we assume that $k>2g-2$. 
The third part of Theorem \ref{thm:Gerbe} is obvious, because the automorphism group of a holomorphic line bundle on a compact Riemann surface is $\bC^{\times}$. The other two parts are in the lemmata below.

\begin{lem}\label{existencelocalsections}
The map $\Phi_g^k: \HolStack_g^k \to \PicStack_g^k$ admits local sections.
\end{lem}

\begin{proof}
Let $X$ be a space, $x \in X$ and $\beta:X \to \PicStack_g^k$ be a map, giving a family of Riemann surfaces $E \to X$ and a section $s: X \to \Pic^k (E/X)$.
We use a classical construction from the geometry of algebraic curves. The fibrewise $k$-fold symmetric product is denoted by $\Sym^k (E/X)\to X$. It is well-known and not hard to see that this is a fibre bundle with smooth complex manifolds as fibres. Recall the classical divisor--line-bundle correspondence \cite[p.\ 129 ff.]{GrifHarr}. This construction yields a fibre-preserving map $\eta:\Sym^k (E/X) \to \Pic^k (E/X)$. 
As long as $k \geq g$, $\eta$ is surjective (the Jacobi inversion theorem; because any line bundle of degree $k$ then has a nontrivial holomorphic section). It is a classical result that for $k> 2g-2$, $\eta$ is a submersion with fibres isomorphic to $\cp^{k-g}$. More precisely, Mattuck \cite{Mattuck} has shown that for an individual Riemann surface $S$, the map $\Sym^k (S) \to \Pic^k (S)$ is a projective bundle (the structural group is $\pgl_{k+1-g}(\bC)$).
It follows that $\eta:\Sym^k (E/X) \to \Pic^k (E/X)$ is a proper surjective submersion. In particular, $\eta$ has local sections and by passing to a smaller $X$ we can pick a lift of $s$; denoted by $t:X \to \Sym^k (E/X)$; this is a section of the bundle $\Sym^k (E/X) \to X$. 

Given a divisor $D$ on a Riemann surface $S$, the classical construction gives an actual line bundle (and not merely an isomorphism class) when one specifies local holomorphic functions on $S$ that define $D$, i.e. have the same zeroes (with mutliplicity). Observe that such local functions can be picked continuously when the divisor varies continuously. This argument shows that after passage to an even smaller $X$, we can find a line bundle $L \to E$ whose restriction to $E_y$ corresponds to the divisor $t(y)$, for all $y \in X$. This $L$ is the desired Poincar\'e line bundle.
\end{proof}

\begin{lem}
Let $\beta:X \to \PicStack_g^k$ and let $L_i \to E$, $i=0,1$, be two Poincar\'e line bundles for $\beta$. Then each $x \in X$ has a neighborhood $U$ such that $(L_0)|_U \cong (L_0)|_U$.
\end{lem}

\begin{proof}
Consider the line bundle $H:=\Hom (L_{0}, L_1) \to E$. It has fibrewise degree $0$, and the hypothesis states that for each $b \in B$, the bundle $H|_{E_b}$ is trivial. Now look at the fibrewise Cauchy--Riemann operator for $H$. As the kernels all have dimension $1$, they form a line bundle on $B$ by basic Fredholm theory. By passing to a neighborhood $U$ of a given $b \in B$, we can pick a nowhere zero section of the kernel bundle. This is nothing else than an isomorphism $L_0|_U \cong L_1|_U$.
\end{proof}

\begin{rem}
One can show a stronger statement. Namely, the fibre bundle $\eta:\Sym^k (E/B) \to \Pic^k (E/B)$ with fibre $\cp^{k-g}$ has structure group $\pgl_{k-g+1}(\bC)$ (if $k>2g-2$).
There exists a cartesian diagram:
\begin{equation*}
\xymatrix{
\HolStack_g^k \ar[r] \ar[d] & \ast \hq \Gl_{k-g+1}(\bC) \ar[d]\\
\PicStack_g^k \ar[r] & \ast \hq \pgl_{k-g+1}(\bC).
}
\end{equation*}
The lower horizontal map classifies the bundle $\eta$ and the upper horizontal map classifies the vector bundle on $\HolStack_g^k$ whose fibre at $(C,L)$ is the space of holomorphic sections of $L$. We do not make use of this result, but leave a remark. The Dixmier--Douady class of the gerbe $\Phi_g^k$ is an element of $H^3 (\PicStack_g^k;\bZ) \cong H^3 (\Pic_g^k;\bZ)$, and under the homotopy equivalence $\Pic_g^k \simeq B \exmcg_{g}(k)$ it corresponds to the generator $\theta \in H^3 (B \exmcg_{g}(k);\bZ)$. Since the order of the gerbe $\ast \hq \Gl_{k-g+1}(\bC)\to \ast \hq \pgl_{k-g+1}(\bC)$ is $k-g+1$, this can be used to give another proof that $\gcd (k-g+1,2g-2) \theta =0$. 
\end{rem}

\begin{proof}[Proof of Proposition \ref{local-quotient-stacks}]
An atlas for $\PicStack_g^k$ is given as follows. Let $\pi:\cE_g \to \cT_g$ be the universal family of Riemann surfaces on Teichm\"uller space. The total space of $\Pic(\cE_g/\cT_g)$ is a complex manifold of dimension $4g-3$. There is a map $\mu:\Pic (\cE_g) \to \PicStack_g$ that fits into a cartesian diagram
\[
\xymatrix{
\Pic (\cE_g) \ar[r] \ar[d] & \PicStack_g \ar[d]\\
\cT_g \ar[r] & \ModStack_g.
}
\]
Therefore $\mu$ is an atlas, and so $\PicStack_g$ is a Deligne--Mumford stack. The corresponding result for $\HolStack_g^k$ follows from this and \cite[Corollary 6.3]{Hein}.
\end{proof}

\subsection{Proof of Theorem \ref{thm:HolCxStacks}}

We will use the basic notions of homotopy theory for stacks as developed in \cite{EbGian} or \cite{Noo}. The homotopy equivalence goes via a zig-zag, and we have to introduce an intermediate stack for this purpose.

\begin{defn}
Let $\CxStack_g$ denote the following stack, defined on the site $\Top$. An object of $\CxStack_g(Y)$ consists of an oriented smooth surface bundle $\pi : E \to Y$ with genus $g$ fibres, and a complex line bundle $L \to E$. An isomorphism $(E, L) \to (E', L')$ is an isomorphism of surface bundles $f : E \to E'$ and an isomorphism $h : L \to f^*L'$ of complex line bundles. We denote by $\CxStack_g^k \subset \CxStack_g$ the open and closed substack consisting of those bundles having fibrewise degree $k$.
\end{defn}

There is a diagram
\begin{equation*}
\begin{gathered}
\xymatrix{
\HolStack_g \ar[r] \ar[d] & \CxStack_g \ar[d] & \modspace{g} \ar[l] \ar[dl] \\
\ModStack_g \ar[r]& \ast \hq \Diff{g}
}
\end{gathered}
\end{equation*}
where the left horizontal maps forget the holomorphic structure and the right horizontal one is given by taking the induced bundle. The vertical morphisms forget the line bundle data. It is clear that the diagram commutes (up to $2$-isomorphism).

\begin{lem}\label{lemmaone}
The map $\modspace{g} \to \CxStack_g $ is a universal weak equivalence.
\end{lem}
\begin{proof}
Since by definition
\[
\modspace{g} = E \Diff{g} \times_{\Diff{g}} \map (\Sigma_g , \CPinf);
\]
we can rewrite the map $\modspace{g}  \to \CxStack_g$ as the composition
\begin{equation}\label{eq:composition}
E \Diff{g} \times_{\Diff{g}} \map (\Sigma_g , \CPinf) \to \map (\Sigma_g , \CPinf)\hq \Diff{g} \to \CxStack_g.
\end{equation}
The first map is a universal weak equivalence by general stack-theoretic principles, namely \cite[Proposition 2.5]{EbGian}. To analyze the second map, consider a space $X$ and a map $X \to \CxStack_g$, representing $L \to E \stackrel{\pi}{\to} X$, where $\pi:E \to X $ is a surface bundle and $L \to E$ a complex line bundle.

Look at the auxiliary space $S$ of pairs $(x,f)$, where $x \in X$ and $f: L|_{\pi^{-1}(x)} \to \pi^{-1}(x) \times \bC^{\infty}$ is a bundle monomorphism. The topology on $S$ is the unique one such that the projection $S \to X$ is a locally trivial fibre bundle and such that the preimage of $x$ has the compact-open topology. The fibre over $x$ is the space of bundle monomorphisms $L|_{\pi^{-1}(x)} \to \pi^{-1}(x) \times \bC^{\infty}$, which is contractible.
It is a routine verification to identify the fibre product $X \times_{\CxStack_g} \map (\Sigma_g ; \CPinf)\hq \Diff{g} $ with $S$. Thus the second map in (\ref{eq:composition}) has local sections and it is a universal weak equivalence. This finishes the proof
\end{proof}

\begin{lem}\label{lemmatwo}
The map $\HolStack_g \to \CxStack_g $ is a universal weak equivalence.
\end{lem}

\begin{proof}
The map $\HolStack_g \to \CxStack_g$ can be factored into three maps
\[
\HolStack_g \stackrel{\phi_1}{\to} \widetilde{\HolStack}_g\stackrel{\phi_2}{\to} \qHolstack_g \stackrel{\phi_3}{\to} \CxStack_g,
\]
each of which is a homotopy equivalence. The stack $\qHolStack_g$ parametrises families of Riemann surfaces together with complex line bundles, and the map $\phi_3$ is the forgetful map. 
It is a universal weak equivalence because the diagram
\begin{equation*}
\xymatrix{
\qHolstack_g \ar[d] \ar[r] & \CxStack_g \ar[d] \\
\ModStack_g \ar[r] & \ast \hq \Diff{g}
}
\end{equation*}
is a fibre square and the bottom map is a universal weak equivalence by Teichm\"uller thory.

Now we define the stack $\widetilde{\HolStack}_g$. A map $X \to \widetilde{\HolStack}_g$ is an element $L \to E \to X$ of $\qHolStack_g (X)$, together with a family of fibrewise differential operators $D: \sect (E; L) \to \sect (E; \Lambda^{0,1}_{v} \otimes L)$ such that 
\begin{equation}\label{symbolequation}
D(fs) = \delbar f \otimes s + f D(s)
\end{equation}
for each $s \in \sect (E; L)$ and $f \in C^{\infty} (E)$, where differentiation is understood to be in the fibrewise sense. The map $\phi_2$ forgets the differential operator. The condition (\ref{symbolequation}) is convex, which implies that each smooth line bundle admits such an operator and that the space of these operators is convex. Therefore $\phi_2$ is a universal weak equivalence. 

The map $\phi_1$ associates to each holomorphic line bundle the Cauchy--Riemann operator on that line bundle. We claim that $\phi_1$ is an isomorphism of stacks. This amounts to showing that a holomorphic structure on a line bundle is determined by its Cauchy--Riemann operator (which is a tautology) and that any family of operators satisfying (\ref{symbolequation}) induces the structure of a holomorphic line bundle on $L$, i.e.\ there exist locally nonzero solutions of $Ds=0$. The argument we give is due to Atiyah and Bott \cite[p.\ 555]{AtiyahBott:Yang-Mills} (they consider the case of higher-dimensional vector bundles, which is more complicated, but for a fixed Riemann surface, which is easier).

Let $(L \to E\stackrel{p}{\to} X ,D)$ be an element of $\widetilde{\HolStack}_g(X)$. For $x \in X$ and $y \in p^{-1}(x)$, we can pick a neighborhood $U$ of $x$ and a map $\alpha:U \times \bD \to E$ over $X$ which is an open embedding, fibrewise holomorphic and satisfies $\alpha(x,0)=y$; furthermore we require $L$ to be trivial over $U \times \bD$. We wish to find a section $s$ of $L$ over $\alpha(U \times \bD)$ that is nowhere zero and satisfies $Ds=0$ over $U \times \bD_{\frac{1}{2}}$. To this end, pick a fibrewise smooth section $s_0$ of $L$ over $\alpha(U\times \bD)$ and look for a function $f$ that satisfies $D(e^f s_0)=0$. So we have to solve the PDE
\[
0= e^{-f} D(e^{f}s_0)=\delbar f \otimes s_0 + Ds_0.
\]
Write $Ds_0=-\beta \otimes s_0$ for a $(0,1)$-form $\beta$; this reduces the problem to the equation
\[
\delbar (f) = \beta.
\]
Since all that matters is a local section on $U \times \bD_{\frac{1}{2}}$, we can multiply $\beta$ with a cut-off function and thus assume that $\beta$ has compact support. Now we pick a fibrewise holomorphic embedding $U \times \bD \to U \times \cp^1$ and a bundle map from the trivial bundle on $U \times \bD$ into the tautological line bundle on $U \times \cp^1$. By Riemann--Roch, the Cauchy--Riemann operator on the tautological line bundle (it has degree $-1$) is invertible. Hence its inverse is continuous as well. Therefore we can find a continuous solution of $\delbar f = \beta$ over $U \times \bD$. The arguments given so far amount to the construction of an inverse map to $\phi_1$ and thus the proof is complete.
\end{proof}

Now we will show how to compare $\PicStack_g^k$ with $B\exmcg_g(k)$. The first observation is that the homotopy type $\Pic_g^k$ is aspherical, as the map $\PicStack_g^k \to \ModStack_g$ is representable and a torus bundle, and the homotopy type of $\ModStack_g$ is $B\Gamma_g$, so aspherical.

\begin{lem}
After taking homotopy types, the map $\Phi_g^k : \HolStack_g^k \to \PicStack_g^k$ is a first Postnikov approximation.
\end{lem}
\begin{proof}
As the homotopy type of $\PicStack_g^k$ is aspherical, it is enough to show that the map $\Phi_g^k$ induces an isomorphism on fundamental groups of homotopy types at all basepoints. By \cite[Example 4.4 and Theorem 5.2]{Noo2}, there is an exact sequence
$$\cdots \lra \pi_n(\CPinf, *) \lra \pi_n(\Hol_g^k, x) \lra \pi_n(\Pic_g^k, \phi_g^k(x)) \lra \pi_{n-1}(\CPinf, *) \lra \cdots$$
on homotopy groups. As $\CPinf$ is simply connected, the claim follows.
\end{proof}

We have the diagram of solid arrows
\[
\xymatrix{
\modspace{g}(k) \ar[r]^{\simeq} \ar[d]^{\Pi}& \mathsf{Cx}_g^k  & \HolStack_g^k \ar[l]_{\simeq} \ar[d]^{\Phi_g^k}\\
B \exmcg_g(k) \ar@{..>}[rr] & & \PicStack_g^k
}
\]
By taking associated homotopy types, and using obstruction theory, there exists a dotted map (on homotopy types) which is a weak equivalence. Composing with $\Pic_g^k \to \PicStack_g^k$, we obtain the map in the statement of


\subsection{Holomorphic versus topological Picard groups}\label{sec:AG}

\begin{thm}\label{comparison:linebundles}
For both stacks $\PicStack_g^k$ and $\HolStack_g^k$, the comparison map $\Pic_{\hol}(-) \to H^2 (- \; ;\bZ)$ is an isomorphism.
\end{thm}

\begin{proof}[Proof of surjectivity]
We know that $H^2 (\HolStack_g^k;\bZ) \cong \bZ^3$ generated by $\lambda$, $\mclass_{0,1}$ and $\zeta$, by Theorem \ref{modspaceproperties} (\ref{it:stability}), Theorem \ref{thm:LowDimCohomologyHol}, and Theorem \ref{thm:HolCxStacks}. The index bundles of the Cauchy--Riemann operators on powers of the tangent bundle and the line bundle are holomorphic vector bundles, and in equation (\ref{grothendieckriemannroch}) we computed their first Chern classes. 

To realise $\lambda$ as a first Chern class, put $r=s=0$. To realise $\lambda + \zeta$, put $(r,s)=(0,1)$. Finally, $s=r=1$ yields $ 13 \lambda+\zeta + \mclass_{0,1}$, which completes the proof for $\mathsf{X}=\HolStack_g^k$. The result for $\PicStack_g^k$ follows from this and Lemma \ref{descend-line-bundles-in-gerbes} below.
\end{proof}

\begin{lem}\label{descend-line-bundles-in-gerbes}
Let $\mathsf{X}$ and $\mathsf{Y}$ be holomorphic local quotient stacks and $\mathsf{X} \to \mathsf{Y}$ a $\bC^{\times}$-gerbe. If the comparison map $\Pic_{\hol}(\mathsf{X})\to \Pic_{\top}(\mathsf{X})$ is surjective, then so is $\Pic_{\hol}(\mathsf{Y})\to \Pic_{\top}(\mathsf{Y})$.
\end{lem}
\begin{proof}
Any line bundle (holomorphic or topological) $L\to X$ has an action of $\bC^{\times}$ coming from the gerbe structure. As usual, $z \in \bC^{\times}$ acts by multiplication with $z^w$ for a uniquely determined $w \in \bZ$, called the \emph{weight} of $L$. 
It is not difficult to see that $w: \Pic_{\top}(\mathsf{X}) \to \bZ$ is a homomorphism and coincides with the edge homomorphism $H^2 (X;\bZ) \to H^2 (\CPinf;\bZ)$ derived from the Leray--Serre spectral sequence. Moreover, a line bundle $L \to \mathsf{X}$ descends to a line bundle on $\mathsf{Y}$ if and only if $w(L)=0$.
Therefore we get a diagram with exact rows:
\begin{equation*}
\xymatrix{
  & \Pic_{\hol}(\mathsf{Y}) \ar[r] \ar[d] & \Pic_{\hol}(\mathsf{X}) \ar[r]^-{w} \ar@{->>}[d] & \bZ \\
0  \ar[r] & \Pic_{\top} (\mathsf{Y})\ar[r]   & \Pic_{\top}(\mathsf{X}) \ar[r]^-{w} & \bZ.
}
\end{equation*}
The bottom sequence is exact at the left by a look at the Leray--Serre spectral sequence. A fragmentary version of the 5-lemma holds in this situation and shows that the left vertical map is onto.
\end{proof}

\begin{proof}[Proof of injectivity]
This begins with a result by Arbarello and Cornalba \cite{AC}. They showed that $\Pic_{\hol}(\ModStack_{g}^{n}) \to H^2 (\ModStack_{g}^{n};\bZ)$ is an isomorphism for all $n \geq 0$. The long exact cohomology sequence of the exponential sequence, using that $H^1 (\ModStack_{g}^{n}; \bZ)=0$, shows $H^1 (\ModStack_{g}^{n}; \mathcal{O})=0$ for all $n \geq 0$. We use the Leray spectral sequence for the map $f: \PicStack_g^k \to \ModStack_g$, which takes the form
$$E_{2}^{p,q} := H^p (\ModStack_g; R^q f_* \cO) \Longrightarrow H^{p+q}(\PicStack_{g}^{k};\cO).$$

As $f$ is a proper map with connected fibres, $f_* \mathcal{O} = \mathcal{O}$, and so $E_2^{1,0}=0$.

Let $\pi : \ModStack_g^1 \to \ModStack_g$ denote the universal curve. There is a map $i : \ModStack_g^1 \to \PicStack_g^k$ over $\ModStack_g$ given by $(\Sigma, p) \mapsto (\Sigma, K + (k - 2g+2)\cdot p)$ where $K$ denotes the canonical class of $\Sigma$. It satisfies $f \circ i = \pi$, and
$$i^* : R^1 f_* \mathcal{O} \lra R^1 \pi_* \mathcal{O}$$
is known to be an isomorphism. We argue that $H^0 (\ModStack_g; R^1 \pi_* \cO)=0$. Otherwise, there would be a nontrivial section of the holomorphic vector bundle $R^1 \pi_* \cO$, producing a nontrivial family of sections of $\Pic^0 (\ModStack_{g}^{1}/\ModStack_g)$, parametrised by $\bC$. This would produce a line bundle on $\ModStack_{g}^{1}$ which has fibrewise degree $0$, but is nontrivial in some fibre, contradicting the result of Arbarello and Cornalba just quoted. Therefore $H^1 (\PicStack_g^k;\mathcal{O})=0$. 

To derive that $H^1 (\HolStack_g^k;\mathcal{O})=0$ as well, consider the Leray--Serre spectral sequence of the gerbe $\HolStack_g^k \to \PicStack_g^k$, which has the form $E_{2}^{p,q}= H^p (\PicStack_g^k; H^q (\ast \hcoker \bC^{\times};\mathcal{O}))$. Since $\bC^{\times}$ is reductive, this $E_2$-term is concentrated in the bottom row and because $H^0 (\ast \hcoker \bC^{\times};\mathcal{O})=\bC$, the proof is complete.
\end{proof}

\subsection{Relation to algebraic geometry}
In the algebraic setting there exists a smooth algebraic stack $\mathcal{P}\mathrm{ic}_g^k$ over the moduli stack of curves $\mathcal{M}_g$, which is an algebraic analogue of our $\HolStack_g^k \to \ModStack_g$. There is a natural copy of the multiplicative group $\bG_m$ inside the endomorphisms of every object of $\mathcal{P}\mathrm{ic}_g^k$, hence it may be rigidified to a stack $\mathcal{P}_g^k$ which is an algebraic analogue of our $\PicStack_g^k$. This new stack is smooth and Deligne--Mumford, and the quotient map is representable. We have not proved, but expect to be true, that the associated analytic stacks to  $\mathcal{P}\mathrm{ic}_g^k$ and $\mathcal{P}_g^k$ are $\HolStack_g^k$ and $\PicStack_g^k$ respectively. In any case, there is certainly an analytic map $(\mathcal{P}_g^k)^{an} \to \PicStack_g^k$, and so we may consider the composition
\begin{equation}\label{eq:AlgAnalytComparison}
\Pic(\mathcal{P}_g^k) \lra \Pic_{\hol}((\mathcal{P}_g^k)^{an}) \lra \Pic_{\hol}(\PicStack_g^k),
\end{equation}
and similarly for $\mathcal{P}\mathrm{ic}_g^k$.

In \cite{MeloViviani}, Melo and Viviani use algebro-geometric methods to study the Picard groups of $\mathcal{P}_g^k$ and $\mathcal{P}\mathrm{ic}_g^k$, and also of their compactifications. Comparing with \cite[Theorem 4.2 and Corollary 4.4]{MeloViviani}, we see---by computation of both sides---that the composition (\ref{eq:AlgAnalytComparison}) and its analogue for $\mathcal{P}\mathrm{ic}_g^k$ are both isomorphisms.

\subsection{A technical result on local quotient stacks}\label{subsct:acyclicity}

\begin{lem}\label{lem:TopPicardGpInj}
Let $\mathsf{X}$ be a local quotient stack.
Then the sheaf $\mathcal{C}_\mathsf{X}$ of continuous complex-valued functions is acyclic, and so the homomorphism $c_1 : \Pic_{\top}(\mathsf{X}) \to H^2(\mathsf{X};\bZ)$ is an isomorphism.
\end{lem}
\begin{proof}
This follows from a twofold application of the descent spectral sequence. If $\mathsf{X}$ and $\mathsf{Y}$ are stacks, $\mathsf{Y} \to \mathsf{X}$ a representable surjective map and $\cF$ a sheaf on $\mathsf{X}$, then there is a spectral sequence 
$$E^{p,q}_{1} = H^q (\mathsf{Y}_p; \cF_p) \Longrightarrow H^{p+q}(\mathsf{X};\cF),$$
where $\mathsf{Y}_p := \mathsf{Y} \times_{\mathsf{X}} \mathsf{Y} \times_{\mathsf{X}} \ldots  \times_{\mathsf{X}}  \mathsf{Y}$ ($p$ factors) and $\cF_p$ is the pullback of $\cF$ to $\mathsf{Y}_p$. 

If $\mathsf{X}=X \hq G$ is a global quotient of a space by a compact Lie group, then the descent spectral sequence is 
$$
E_2^{p,q} = H^p_{cts}(G; H^q(X, \mathcal{C}_X)) \Longrightarrow H^{p+q}(\mathsf{X};\mathcal{C}_\mathsf{X});
$$
this is zero if $q>0$ since the sheaf $\mathcal{C}_X$ is fine. Thus the spectral sequence collapses to the continuous group cohomology $H^{*}_{cts}(G;\map(X, \bC))$, which vanishes in positive degrees as the group $G$ is compact and the coefficient module is a locally convex topological vector space with a continuous $G$-action. For details of that argument, consult \cite[Proposition 6.3]{AtiyahSegal}.

If $\mathsf{X}$ is merely a local quotient stack, there is an open cover by substacks $\mathsf{X}_i$ each of which is a global quotient stack. Since fibre products of global quotients are global quotient stacks, the descent spectral sequence for the map $\coprod_{i} \mathsf{X}_i\to \mathsf{X}$ has $E^{p,q}_{1}=0$ for $q>0$. The $E^{p,0}$-line is just the {\v C}ech complex for the open cover. Because there exist partitions of unity in this situation, see \cite[Appendix A]{EbGian}, this complex is acyclic.
\end{proof}

\end{document}